\documentclass[smallextended,numbook,runningheads]{svjour3}
\journalname{bit}
\usepackage{graphicx,enumitem,dsfont}
\usepackage[usenames]{color}
\usepackage{dcolumn}
\usepackage[colorlinks]{hyperref}

\usepackage{mathptmx}

\usepackage{amsmath}
\usepackage{amssymb}
\usepackage{amsfonts}

\usepackage{subfigure}

\newcommand{\R}{\mathbb{R}}
\newcommand{\dd}{\mathbf{d}}

\newcommand{\U}{{\mathbf{u}}}
\newcommand{\B}{{\mathbf{B}}}
\renewcommand{\P}{{\mathbf{P}}}
\newcommand{\Q}{{\mathbf{Q}}}
\newcommand{\X}{{\mathbf{x}}}
\newcommand{\E}{{\mathbf{E}}}
\newcommand{\s}{{\mathbf{S}}}
\newcommand{\g}{{\mathbf{g}}}

\newcommand{\Div}{\mathrm{div}}
\newcommand{\Curl}{\mathrm{curl}}

\newcommand{\Comment}[1]{}

\renewcommand{\i}{\ifmmode\mathit{\mathchar"7010 }\else\char"10 \fi}
\renewcommand{\j}{\ifmmode\mathit{\mathchar"7011 }\else\char"11 \fi}

\newcommand{\seq}[1]{\left\{#1\right\}}

\newcommand{\test}{\varphi}
\newcommand{\Dx}{\Delta x}
\newcommand{\Dy}{\Delta y}

\newcommand{\norm}[1]{\left\|#1\right\|}
\newcommand{\abs}[1]{\left|#1\right|}

\newcommand{\charf}{\mathbf{1}}
\newcommand{\Hdiv}{H^{\mathrm{div}}}
\newcommand{\mx}[2]{\left(#1 \vee #2\right)}
\newcommand{\mi}[2]{\left(#1 \wedge #2\right)}
\newcommand{\e}[1]{\ensuremath{\cdot 10^{#1}}}

\newcounter{asnr}

{\ifnum\value{asnr}=0 \stepcounter{asnr} 
  \begin{enumerate}[label=\textbf{A}.\arabic{enumi}]
    \else
    \begin{enumerate}[label=\textbf{A}.\arabic{enumi},resume] \fi}
{\end{enumerate}}

\begin{document}

\title{Higher order finite difference schemes for the Magnetic Induction equations}

\titlerunning{SBP-SAT schemes for the Magnetic induction equations}

\author{U. Koley \and S. Mishra \and N. H. Risebro \and {M. Sv\"{a}rd}}

\institute{U. Koley \and S. Mishra \and N. H. Risebro \and M. Sv\"{a}rd
\at
Centre of Mathematics for Applications (CMA) \\
University of Oslo \\
P.O. Box 1053, Blindern \\
N--0316 Oslo, Norway
\and 
U. Koley 
\at 
\email{ujjwalk@cma.uio.no}
\and 
S. Mishra
\at 
\email{siddharm@cma.uio.no}
\and
N. H. Risebro
\at
\email{nilshr@math.uio.no}
\and
M. Sv\"{a}rd
\at
SINTEF ICT \\
Pb. 124 Blindern \\
N--0314 Oslo, Norway \\
\email{Magnus.Svard@sintef.no}
}









\date{\today}

\maketitle


\begin{abstract}
  We describe high order accurate and stable finite difference schemes
  for the initial-boundary value problem associated with the magnetic
  induction equations. These equations model the evolution of a
  magnetic field due to a given velocity field. The finite difference
  schemes are based on Summation by Parts (SBP) operators for spatial
  derivatives and a Simultaneous Approximation Term (SAT) technique
  for imposing boundary conditions. We present various numerical
  experiments that demonstrate both the stability as well as high
  order of accuracy of the schemes.


\end{abstract}

\section{Introduction}
\label{sec:intro}
The magnetic induction equations are a special form of the Maxwell's
equations that describe the evolution of the magnetic field under the
influence of a given velocity field. These equations arise in a wide
variety of applications in plasma physics, astrophysics and electrical
engineering. One important application are the equations of
magneto-hydro dynamics (MHD). These
equations combine the Euler equations of gas dynamics with the
magnetic induction equations. Our goal in this paper is to describe
stable and high-order accurate numerical schemes for the magnetic
induction equations.

We start with a brief description of how the equations are derived.
Let the magnetic field and given velocity field be denoted by $\B$ and
$\U$ respectively. Faraday's law for the magnetic flux across a
surface ${\bf S}$ bounded by a curve $\partial {\bf S}$ is given by (see \cite{Pow2}),
\begin{equation}
\frac {d}{dt}\int\limits_{S} \B \cdot \dd\s = \oint\limits_{\partial S}\E \cdot dl . \nonumber
\end{equation}
Using the Stokes theorem and the fact that the electric field, $\E$,
in a co-moving frame is zero and the magnetic resistivity is zero, Faraday's law takes the form,
\begin{equation}
  \label{eq:induc}
  \frac {\partial \B}{\partial t} + {\rm div} ( \U\otimes \B -
  \B\otimes \U) = - \U {\rm div}(\B). 
\end{equation}
Using simple vector identities, (\ref{eq:induc}) can be rewritten as,
\begin{equation}
  \label{eq:induction}
 \partial_t \B + \Curl(\B\times\U) =  - \U \Div(\B). 
\end{equation}
Magnetic monopoles have never been observed in nature. As a
consequence, the magnetic field is always assumed to be divergence
free, i.e., $\Div(\B)=0$. Hence, it is common to set the right-hand
side of (\ref{eq:induction}) to zero and couple the induction equation
with the divergence constraint in order to obtain
\begin{equation}
  \label{eq:Maxwell3D}
  \begin{aligned}
    \partial_t \B + \Curl(\B \times \U) &= 0,\\
    \Div(\B) &= 0, \,
\B(x,0) = \B_0(x).
  \end{aligned}
\end{equation}
This form (\ref{eq:Maxwell3D}) is commonly used in the literature as
the appropriate form of the magnetic induction equations to study and
discretize. It is easy to see that (\ref{eq:Maxwell3D}) is hyperbolic
but not strictly hyperbolic. An important tool in the analysis of
hyperbolic system of equations is the derivation of energy estimates.
The usual procedure in deriving energy estimates consists of
symmetrizing the hyperbolic system. It is not possible to symmetrize
(\ref{eq:Maxwell3D}) without explicitly using the divergence
constraint.  Hence, it is difficult to obtain energy stability
starting from (\ref{eq:Maxwell3D}).

On the other hand, we can use the following vector identity
\begin{align*}
  \Curl(\B\times \U)&=\B\Div \U - \U\Div (\B) +
  \left(\U\cdot\nabla\right) \B -
  \left(\B\cdot\nabla\right) \U\\
  &=\left(u^1\B\right)_x+\left(u^2\B\right)_y+\left(u^3\B\right)_z -
  \U\Div(\B) - (\B\cdot\nabla)\U,
\end{align*}
and rewrite (\ref{eq:induc}) in the non-conservative symmetric form,
\begin{equation}
  \begin{aligned}
    \partial_t\B + \left(\U\cdot\nabla\right)\B&=
    -\B(\Div\U)+(\B\cdot\nabla)\U\\
    &=M(D\U)\B,
  \end{aligned}
  \label{eq:induc1}  
\end{equation}
where the $D\U$ denotes the gradient of $\U$ and the matrix $M(D\U)$ is given by
$$
M(D\U)=
 \begin{pmatrix}
   -\partial_y u^2 - \partial_z u^3 & \partial_y u^1 & \partial_z
   u^1\\
   \partial_x u^2 & -\partial_x u^1-\partial_z u^3 &+\partial_z u^2\\
   \partial_x u^3 & \partial_y u^3 & -\partial_x u^1 - \partial_y u^2
  \end{pmatrix}.
  $$
  Introducing the matrix,
  $$
  C=-\begin{pmatrix}
    \partial_x u^1 & \partial_y u^1 & \partial_z u^1 \\
    \partial_x u^2 & \partial_y u^2 & \partial_z u^2 \\
    \partial_x u^3 & \partial_y u^3 & \partial_z u^1
\end{pmatrix},
$$
(\ref{eq:induc}) can also be written in the following
``conservative'' symmetric form,
\begin{equation}
  \label{eq:induc2}
  \partial_t\B + \partial_x\left(A^1\B\right)
  +\partial_y\left(A^2\B\right)+\partial_z\left(A^3\B\right)+C\B=0,
\end{equation}
where $A^i=u^iI$ for $i=1,2,3$. Note that the symmetrized matrices in
(\ref{eq:induc2}) are diagonal and that the only coupling in the
equations is through the lower order terms. These symmetrized
forms are in the same spirit as the non-linear symmetrized forms of
MHD equations introduced in \cite{God1}.

Furthermore, by taking divergence on both sides of
(\ref{eq:induction}) we get
\begin{equation}
  \label{eq:divt}
  (\Div(\B))_t + \Div\left(\U\Div(\B)\right) = 0.
\end{equation}
Hence, if $\Div(\B_0(\X))=0$, also $\Div(\B(\X,t))=0$ for $t>0$. This
implies that all the above forms (\ref{eq:induc2}), and
(\ref{eq:Maxwell3D}) are equivalent (at least for smooth solutions).
Introducing the space $\Hdiv$ as
$$
\Hdiv(\R^3)=\seq{\mathbf{w}:\R^3\to \R^3 \; \bigm|\;
  \abs{\mathbf{w}}\in L^2(\R^3),\;\; \Div(\mathbf{w})\in L^2(\R^3)},
$$
we have the following theorem:
\begin{theorem}
  Assume that the velocity field $\U$ is sufficiently smooth, and that
  $\B_0\in \Hdiv(\R^3)$. Then there exists a unique weak solution $\B
  \in C([0,T];\Hdiv(\R^3))$ of (\ref{eq:induc2}).  The solution $\B$
  satisfies the energy estimate,
  \begin{equation*}
    \norm{\B(\cdot,T)}_{\Hdiv(\R^3)} \le C_T\norm{\B_0}_{\Hdiv(\R^3)}
  \end{equation*}
  The constant $C_T$ depends only on the final time $T$. Furthermore, if $\Div(\B_0)=0$, then the physical form
  (\ref{eq:induc}) and the symmetric form (\ref{eq:induc2}) are
  equivalent to the constrained form (\ref{eq:Maxwell3D}), i.e., $\B$
  is also the unique weak solution of (\ref{eq:Maxwell3D}).
\end{theorem}
The proof of the above theorem uses the energy estimate and we will
provide a sketch of the proof for the two-dimensional version of the
equations together with boundary conditions later in this paper.

Even though the magnetic induction equations are linear, the presence
of variable coefficients and lower order terms means that general closed form
solutions are not available.  Hence, one has to design suitable
numerical schemes for these equations. Furthermore, since these
equations appear as a sub-model in the MHD equations, the design of stable
and high-order accurate numerical schemes for the induction equations
can lead to the design of robust schemes for the non-linear MHD
equations.

Most of the attention in the literature has been focused on the
constrained form (\ref{eq:Maxwell3D}). The key issue in the design of
a suitable numerical scheme to approximate (\ref{eq:Maxwell3D}) has
been the treatment of the divergence constraint. A widely used approach has
been to employ projection methods based on a Hodge
decomposition of the magnetic field. A base (finite difference or
finite volume) scheme is used to evolve the magnetic field. The
evolved field, which need not be divergence free, is then corrected
for divergence errors by solving an elliptic equation (see \cite{BlBr1}). The resulting method is computationally expensive, as the elliptic equation has to be solved at every time step.

Another common approach is to discretize (\ref{eq:Maxwell3D}) such
that some particular form of discrete divergence is preserved at each
time step (see \cite{TF1}). This approach is equivalent to staggering the velocity and magnetic fields in each direction (see \cite{DW1,BS1,RMJF1,EH1} and a detailed comparison in \cite{Toth1}). Some of these schemes are proved
to be von Neumann stable in the special case of constant
velocity fields. No stability analysis is available either in the
case of variable velocity fields or for problems with boundary
conditions. These schemes also involve wider stencils than what is
required for a standard finite difference scheme.

Despite all the attempts at finding a suitable discretization of
(\ref{eq:Maxwell3D}) and preserving a special form of discrete
divergence, it is not clear as to whether such an approach is
appropriate. Furthermore, there are
many different choices for the discrete divergence operator and
preserving some form of discrete divergence exactly does not lead to
preservation or even keeping divergence errors small for a
different form. The main aim should be to design a stable scheme to
approximate magnetic fields and it is not clear whether preserving
divergence in a particular discrete form helps. One reason for the difficulties in proving stability
of discretizations for (\ref{eq:Maxwell3D}) with general velocity
fields may lie in the very form of these equations. As remarked
earlier, (\ref{eq:Maxwell3D}) are not symmetrizable directly and thus
one cannot obtain energy estimates in this form. This remains true for
discretizations of (\ref{eq:Maxwell3D}).

A different approach consisting of discretizing the physical form
(\ref{eq:induc}) was proposed in \cite{Pow1} for the non-linear MHD
equations. Adapting this to (\ref{eq:induc}) implies using a standard
upwind scheme for the convection part and a centered discretization of
the source terms. From (\ref{eq:divt}), one can expect that divergence
errors will be transported out of the domain for transparent boundary conditions. This approach does not
imply stability either and can lead to oscillations as reported in
\cite{fkrsid1}. A discontinuous Galerkin based discretization of the
symmetric form (\ref{eq:induc2}) was proposed in \cite{BeKr1}.

In a recent paper \cite{fkrsid1}, the authors discretized the
symmetric form (\ref{eq:induc1}) by using a first order accurate upwind finite difference
scheme. The resulting scheme also implied an upwind discretization of
the convection term in (\ref{eq:induc}) with an upwind discretization
of the source term. This scheme was shown to be energy stable even
with variable velocity fields and to be
TVD for constant velocity fields.

Furthermore, boundary conditions were not considered either
in \cite{fkrsid1} or any of the aforementioned papers. High-order
accurate schemes will lead to much better resolution of interesting
solution features and a stable discretization of the boundary
conditions (while still preserving high order of accuracy) is desirable. 

Our aim in this paper is to design stable and high-order
accurate schemes for initial-boundary value problems corresponding to
the magnetic induction equations by discretizing the non-conservative symmetric form (\ref{eq:induc1}). The spatial derivatives are
approximated by second and fourth-order SBP (Summation-By-Parts)
operators. The boundary conditions are weakly imposed by using a SAT
(Simultaneous Approximation Term) and time integration is performed by
standard Runge-Kutta schemes. The SBP-SAT framework has been used to obtain stable and accurate high order schemes for a wide variety of hyperbolic problems in recent years. See
\cite{SN1} and the references therein for more details.

The SBP-SAT schemes use centered finite difference
stencils in the interior, which lead to oscillations in the vicinity
of discontinuities. We apply
well-known SBP-SAT compatible numerical diffusion operators
in case of discontinuous data.

The rest of this paper is organized as follows: In
Section~\ref{sec:cp}, we state the energy estimate for the
initial-boundary value problem corresponding to (\ref{eq:induc1}) in
order to motivate the proof of stability for the scheme. In
Section~\ref{sec:scheme}, we present the SBP-SAT scheme and show stability. Numerical experiments are presented in
Section~\ref{sec:numex} and conclusions are drawn in
Section~\ref{sec:conc}.

\section{The Continuous problem}
\label{sec:cp}
For ease of notation, we shall restrict ourselves to two spatial
dimensions in the remainder of this paper. Extending the results to
three dimensions is straightforward.

In two dimensions, the non-conservative symmetric form
(\ref{eq:induc1}) reads
\begin{equation}
  \label{eq:main}
  \begin{aligned}
    \B_t + \Lambda_1\B_x + \Lambda_2\B_y - C\B &= 0,
    \end{aligned}
\end{equation}
where
\begin{align*}
  \Lambda_1&=
  \begin{pmatrix}
    u^1 & 0 \\
    0 & u^1
  \end{pmatrix},\quad
  \Lambda_2=
  \begin{pmatrix}
    u^2 & 0 \\
    0 & u^2
  \end{pmatrix},\quad
  C=\begin{pmatrix}
    -\partial_y u^2 & \partial_y u^1 \\
    \partial_x u^2 & -\partial_x u^1
  \end{pmatrix},
\end{align*}
with $\B = \left( B^1, B^2\right)^T$ and $\U = \left(u^1,u^2\right)^T$
denoting the magnetic and velocity fields respectively. In component
form, (\ref{eq:main}) becomes
\begin{equation}
  \label{eq:NC2D}
  \begin{aligned}
    (B^1)_t + u^1(B^1)_x + u^2(B^1)_y &= -(u^2)_y B^1 + (u^1)_y B^2 \\
    (B^2)_t + u^1(B^2)_x + u^2(B^2)_y &= (u^2)_x B^1 - (u^1)_x B^2.
  \end{aligned}
\end{equation}
To begin with, we shall consider (\ref{eq:main}) in the domain
$(x,y)\in\Omega=[0,1]^2$.

We augment (\ref{eq:main}) with initial conditions,
\begin{equation}
  \label{eq:initial}
  \begin{aligned}
    \B(\X,0) = \B_0(\X) \quad \X\in \Omega,
 \end{aligned}
\end{equation}
and Dirichlet boundary conditions,
\begin{equation}
  \label{eq:boundary}
  \begin{aligned}
    \charf_{\seq{u^1(0,y,t)>0}}\Bigl(\B(0,y,t) &= \g(0,y,t)\Bigr),
    \quad
    \charf_{\seq{u^1(1,y,t)<0}} \Bigl(\B(1,y,t) = \g(1,y,t)\Bigr),\\
    \charf_{\seq{u^2(x,0,t)>0}}\Bigl(\B(x,0,t) &=\g(x,0,t)\Bigr),
    \quad \charf_{\seq{u^2(x,1,t)<0}}\Bigl(\B(x,1,t) = \g(x,1,t)\Bigr)
  \end{aligned}
\end{equation}
where $\charf_{A}$ denotes the characteristic function of the set $A$.
Note that we only impose boundary conditions on the set where the
characteristics are entering the domain.

\begin{definition}
Weak solution: A function $\B:\Omega\to \R^2$ such that $\B\in
C([0,T];H^1(\Omega))$ is defined as a weak solution of (\ref{eq:main}) with initial data (\ref{eq:initial}) and boundary data (\ref{eq:boundary}) if it satisfies the weak formulation of
(\ref{eq:main}) in $\Omega$, i.e.,
\begin{multline}
  \label{eq:weakformulation}
  \int_0^T\int_\Omega \B \left(\test_t + \left(\Lambda_1\test\right)_x
    + \left(\Lambda_2\test\right)_y - C\test \right)\, dxdy dt +
  \int_\Omega \B_0 \test(x,y,0)\,dxdy \\ -\int_0^T \int_0^1u^1
  \left(Tr\B\right) \test(x,y,t)\bigm|^{x=1}_{x=0}\,dydt -\int_0^T
  \int_0^1 u^2 \left(Tr\B\right) \test(x,y,t)\bigm|^{y=1}_{y=0}\,dx dt
  = 0,
\end{multline}
for all test functions $\test\in C^\infty_0(\Omega\times [0,T))$. By
$Tr\B$ we mean the $H^1$ trace of $\B$ at the boundary.  The boundary
conditions (\ref{eq:boundary}) are taken in the sense of $H^1$ traces.
\end{definition}
We shall always assume that the initial and boundary data satisfy the
compatibility conditions, i.e., specific criteria that guarantee smoothness of the solution, see \cite{GustafssonKreissOliger}.
\begin{theorem}
  \label{theo:cp}
  Assume that $\B_0\in H^1(\Omega)$, that $\g\in
  H^1(\partial\Omega\times [0,T])$ for $T>0$ and that $u^1$ and $u^2$
  are in $H^2(\Omega\times[0,T])$. Then there exists a function $\B \in
  C([0,T],L^2(\Omega))\cap L^\infty([0,T];H^1(\Omega))$ which is the
  unique weak solution of (\ref{eq:main}) with the initial and
  boundary conditions (\ref{eq:initial}) and (\ref{eq:boundary}).

  Furthermore, it satisfies the following stability estimate
  \begin{equation}
    \label{eq:enest}
    \norm{\B(\cdot,t)}_{H^1(\Omega)}^2\le e^{\alpha t}\left(\norm
    {\B_0}_{H^1(\Omega)}^{2} + \norm{\g}_{H^1(\partial\Omega\times(0,t))}\right).
  \end{equation}
  where $\alpha$ is a positive constant.
\end{theorem}
\begin{proof}
  The proof of this theorem is standard. Assume first that $\g$,
  $\B_0$ and $\U$ are in $C^\infty$. Since the compatibility
  conditions are satisfied, a unique solution exists by the method of
  characteristics. Let $\mx{a}{0}=\max\seq{a,0}$ and
  $\mi{a}{0}=\min\seq{a,0}$. Multiplying the equation by $\B$ and
  integrating over $\Omega$ yields
  \begin{align*}
&\qquad    \frac{d}{dt}\int_\Omega \B^T\B \,dxdy \\
&= \int_\Omega
    \B\left(2C+\Div(\U)\right)\B\,dxdy - \int_0^1 u^1 Tr(\B^T\B)
    \bigm|^{x=1}_{x=0}\,dy +
    \int_0^1 u^2 Tr(\B^T\B) \bigm|^{y=1}_{y=0}\,dx \\
    &\le c \int_\Omega \B^T\B \,dxdy \\
    &\qquad + \int_0^1 \mx{u^1(0,y,t)}{0}
    \left(Tr(\B^T\B)\right)\,dy -
    \int_0^1 \mi{u^1(1,y,t)}{0} \left(Tr(\B^T\B)\right)\,dy\\
    &\qquad + \int_0^1 \mx{u^2(x,0,t)}{0} \left(Tr(\B^T\B)\right)\,dx -
    \int_0^1 \mi{u^2(x,1,t)}{0} \left(Tr(\B^T\B)\right)\,dx\\
    &\le c \biggl(\int_\Omega (\B^T\B) \,dxdy + \int_{\partial\Omega}
    \g^2\,ds \biggr)
  \end{align*}
  for some constant $c$ depending on $\U$ and its first derivatives.
  Via the Gr\"onwall inequality we get the bound
  $$
  \norm{\B(\cdot,t)}_{L^2(\Omega)}^2 \le e^{ct}
  \left(\norm{\B_0}_{L^2(\Omega)}^2 + \int_0^T \int_{\partial\Omega}
    \g^2\,ds\,dt \right).
  $$
  Set $\P=\B_x$ and $\Q=\B_y$, applying $\partial_x$ to
  (\ref{eq:main}) yields
  \begin{equation}
    \P_t + u^1 \P_x + u^2 \P_y = u^1_x \P + u^2_x \Q + C\P + C_x
    \B.\label{eq:Peq}
  \end{equation}
  Furthermore, $P(x,y,0)=\partial_x\B_0(x,y)$ and at those parts of
  $\partial\Omega$ where we impose boundary data
  $$
  \begin{aligned}
    u^1 \P &= C\g - \g_t - u^2 \g_y \quad\text{on $x=0$ and $x=1$,}\\
    u^2 \Q &= C\g - \g_t - u^2 \g_x \quad\text{on $y=0$ and $y=1$.}
  \end{aligned}
  $$
  We shall also be needing $\P$ on $y=0$ and $1$ and $\Q$ on $x=0$
  and $1$. These are given by $\g_x$ and $\g_y$ respectively.
  
  Multiplying (\ref{eq:Peq}) with $2\P^T$ and rearranging yields
  $$
  \P^2_t + \left(u^1\P^2\right)_x + \left(u^2\P^2\right)_y = -u^1_x
  \P^2 -2u^2_x \P^T\Q +2\P^TC\P + 2\P^TC_x\B.
  $$
  We also have an analogous equation for $\Q^2$;
  $$
  \Q^2_t + \left(u^1\Q^2\right)_x + \left(u^2\Q^2\right)_y = -u^2_y
  \Q^2 -2u^1_y \P^T\Q +2\P^TC\Q + 2\Q^TC_y\B.
  $$
  Adding these two equations we find
  \begin{equation}
    \label{eq:PplussQ}
    \left(\P^2+\Q^2\right)_t +
    \left(u^1\left(\P^2+\Q^2\right)\right)_x + 
    \left(u^2\left(\P^2+\Q^2\right)\right)_y = R,
  \end{equation}
  where by H\"older's inequality $R$ has the bound
  $$
  \int_{\Omega}R(x,y,t)\,dxdy \le c\left(
    \norm{\B(\cdot,t)}_{L^2(\Omega)}^2 +
    \norm{\P(\cdot,t)}_{L^2(\Omega)}^2 +
    \norm{\Q(\cdot,t)}_{L^2(\Omega)}^2\right),
  $$
  where the constant $c$ depends on $\U$ and its derivatives.  By
  reasoning as we did with $\B$, we can then get the bound
  \begin{multline*}
    \frac{d}{dt}\left(\norm{\P}_{L^2(\Omega)}^2+\norm{\Q}_{L^2(\Omega)}^2
    \right)\le  \\
    c\left(\norm{\B(\cdot,t)}_{L^2(\Omega)}^2+
      \norm{\P}_{L^2(\Omega)}^2
      +\norm{\Q}_{L^2(\Omega)}^2+\int_{\partial\Omega} \g_t^2 +
      \g_x^2+ \g_y^2\,ds\right).
  \end{multline*}
  Via Gr\"onwall's inequality and the bounds on $\norm{\B}_{L^2}$ we
  find
  $$
  \norm{\P(\cdot,t)}_{L^2(\Omega)}^2
  +\norm{\Q(\cdot,t)}_{L^2(\Omega)}^2\le \mathrm{Const.},
  $$
  where the constant depends on the $H^1(\Omega)$ norm of $\B_0$
  and $\g$ and on $\U$ and its derivatives. This means that we have an
  energy estimate
  $$
  \norm{\B(\cdot,t)}_{H^1(\Omega)} \le C_t
  \left(\norm{\B_0}_{H^1(\Omega)}+ \norm{\g}_{H^1(\partial\Omega\times
      (0,t))}\right),
  $$
  where $C_t$ is a finite constant depending on $t$, $\U$ and its
  derivatives.
  
  Then, for a general initial data, velocity fields and boundary conditions, we can use a standard approximation argument (\cite{KL1})  and the above estimate to pass to the limit and prove the existence and uniqueness of weak solutions. 
\end{proof}
\begin{remark}
The above theorem has been proved in the unit square. It can be easily extended to domains with smooth (i.e., $C^1$ boundaries) by using cut-off functions and mappings between the domain and the upper-half space. See \cite{R} and other references therein for details.
\end{remark} 
\section{Semi-discrete Schemes}
\label{sec:scheme}
As stated before, we will approximate (\ref{eq:main}) with SBP-SAT
finite difference schemes. We start by defining a SBP operator
approximating the first derivative of a continuous function $w(x)$ in
one space dimension. Let $\seq{x_i}_{i=0}^n$ be equidistant points in
$[0,1]$ such that $x_i=ih$ where $h=1/n$.  We organize the values of
$w$ at $\seq{x_i}$ in a vector $w^T(t)=(w_0,...,w_n)$ where
$w_i=w(x_i)$.  Then , we define,
\begin{definition}
  A difference approximation (given by a $(n+1)\times (n+1)$ matrix $D$) for
  the first derivative is called a \emph{Summation-By-Parts} (SBP)
  operator if $D=P^{-1}Q$ for $n\times n$ matrices $P$ and $Q$, where
  $P > 0$, $P=P^T$ and $Q+Q^T=\mathcal{B}=\mathrm{diag}(-1,0,0,\ldots,0,0,1)$.
  
  Moreover, $P$ must define  a scalar product $(w,v)=w^TPv$ for which
  the corresponding norm, $\norm{w}^2_P=(w,w)$, is  equivalent to the
  standard $l^2$-norm, $\norm{w}_2^2=h\sum_{i=1}^n w_i^2$.
\end{definition}
SBP operators of different orders of accuracy are presented in several
papers, see the references in e.g.~\cite{SN1}. To discretize (\ref{eq:main}),
we introduce equidistant meshes in the $x$ and $y$ directions with
$N$ and $M$ mesh points and $\Dx=1/N$ and $\Dy=1/M$. The discrete
solution consists of a column vector of length $2(N+1)(M+1)$
denoted  $V=(V^1,V^2)^T$,  where $V^\ell$ is a vector of length
$(N+1)(M+1)$ ordered as 
$$
V^\ell=\left(
  V^\ell_{0,0},V^\ell_{0,1},\ldots,V^\ell_{0,N},
  V^\ell_{1,0},\ldots,\ldots,V^\ell_{N,M}\right).
$$
and $V^\ell_{i,j}$ is the discrete solution at $(x_i,y_j)$ for
$\ell=1,\,2$. 
We will use the
norm
$$
\norm{w}_{P}^2=w^T(P_x\otimes P_y)w
$$
where we have introduced the Kronecker product, which is defined as
follows:

 Let $A$ and $C$ be $n\times n$ matrices and $B$ and $D$ be
$m\times m$ matrices. Then $A\otimes B$ is the $nm\times nm $ matrix 
$$
(A\otimes B)=
\begin{pmatrix} a_{11}B &\hdots  & a_{1n}B\\
  \vdots & \ddots & \vdots \\
  a_{n1}B &\hdots  & a_{nn}B
\end{pmatrix}.
$$
Furthermore, the following rules hold;
$(A\otimes B)(C\otimes D)=(AC\otimes BD)$, $(A\otimes B)+(C\otimes
D)=(A+C)\otimes(B+D)$ and $(A\otimes B)^T=(A^T\otimes B^T)$.

To define discrete boundary conditions, we need some further
notation. For real numbers $\sigma_i$, introduce the $2\times 2$ matrices
\begin{align*}
  \Sigma_{0,y}=\sigma_1 I_2 ,\quad 
  \Sigma_{N,y}= \sigma_2 I_2,\quad 
  \Sigma_{x,0}=\sigma_3 I_2,\quad 
  \Sigma_{x,N}= \sigma_4 I_2,
\end{align*}
where the $I_2$ is the $2\times 2$ identity matrix and the numbers $\sigma_i$ are to be determined later. We
also need $(M+1)\times(M+1)$ matrices $F_{0,y}$ and $F_{N,y}$
\begin{align*}
  F_{0,y}&=\begin{pmatrix}
    1 & 0 &\cdot &\cdot &\cdot & 0 \\
    0 & \cdot &\cdot &\cdot &\cdot & 0 \\
    \cdot &\cdot &\cdot &\cdot &\cdot &\cdot \\
    1 & 0 &\cdot &\cdot &\cdot & 0 \\
    0 & \cdot &\cdot &\cdot &\cdot & 0 \\
    \cdot &\cdot &\cdot &\cdot &\cdot &\cdot 
  \end{pmatrix},\quad
  F_{N,y}=\begin{pmatrix}
    0 & \cdot &\cdot &\cdot &\cdot &1 \\
    0 & \cdot &\cdot &\cdot &\cdot & 0 \\
    \cdot & \cdot &\cdot &\cdot &\cdot & \cdot \\
    0 & \cdot &\cdot &\cdot &\cdot &1 \\
    0 & \cdot &\cdot &\cdot &\cdot & 0 \\
    \cdot & \cdot &\cdot &\cdot &\cdot & \cdot
  \end{pmatrix},\\
  \intertext{and $(N+1)\times(N+1)$ matrices $F_{x,0}$ and $F_{x,M}$,}
  F_{x,0}&=\begin{pmatrix}
    1 & 0 &\cdot &\cdot &\cdot & 0 \\
    0 & \cdot &\cdot &\cdot &\cdot & 0 \\
    \cdot &\cdot &\cdot &\cdot &\cdot &\cdot \\
    1 & 0 &\cdot &\cdot &\cdot & 0 \\
    0 & \cdot &\cdot &\cdot &\cdot & 0 \\
    \cdot &\cdot &\cdot &\cdot &\cdot &\cdot
  \end{pmatrix},\quad
  F_{x,M}=\begin{pmatrix}
    0 & \cdot &\cdot &\cdot &\cdot &1 \\
    0 & \cdot &\cdot &\cdot &\cdot & 0 \\
    \cdot & \cdot &\cdot &\cdot &\cdot & \cdot \\
    0 & \cdot &\cdot &\cdot &\cdot &1 \\
    0 & \cdot &\cdot &\cdot &\cdot & 0 \\
    \cdot & \cdot &\cdot &\cdot &\cdot & \cdot
  \end{pmatrix}.
\end{align*}
Next set
\begin{align*}
  E_{0,y}=I_x\otimes F_{0,y},\quad E_{N,y}= I_x\otimes F_{N,y},\quad  E_{x0}=
  F_{x,0}\otimes I_y,\quad\text{and}\quad E_{x,M}= F_{x,M}\otimes I_y
\end{align*}
where $I_x$ and $I_y$ are $(N+1)\times(N+1)$ and $(M+1) \times (M+1)$ identity
matrices respectively and define
\begin{align*}
  \Lambda_{x} &= I_2\otimes
  \mathrm{diag}
  \left(  
    u^1_{0,0},u^1_{0,1},\ldots,u^1_{0,N},
    u^1_{1,0},\ldots,\ldots,u^1_{N,M}\right)\\
  \Lambda_{y} &= I_2\otimes
  \mathrm{diag}
  \left(  
    u^2_{0,0},u^2_{0,1},\ldots,u^2_{0,N},
    u^2_{1,0},\ldots,\ldots,u^2_{N,M}\right).
\end{align*} 
Define the matrix $C$ by,
\begin{align*}
  C=\begin{pmatrix}
    -\left(I_x\otimes\left(P_y^{-1} Q_y\right)\right) u^2 
    &\left(I_x\otimes\left(P_y^{-1} Q_y\right)\right) u^1 \\
    \left(\left(P_x^{-1} Q_x\right)\otimes I_y\right) u^2 
    &-\left(\left(P_x^{-1} Q_x\right)\otimes I_y\right)u^1
  \end{pmatrix}.
\end{align*}
Let $g$ be a column vector of the same length as $V$, where we
store the boundary values at the appropriate places. Then we can
describe the SBP-SAT scheme as
\begin{equation}
  \label{eq:discrete}
  \begin{aligned}
    \partial_t &V + 
    \Lambda_{x}\left(I_2\otimes\left(P_x^{-1}Q_x\right)\otimes
      I_y\right) V 
    + 
    \Lambda_{y}\left(I_2\otimes I_x\otimes\left(P_y^{-1}
        Q_y\right)\right)V 
    +CV\\
    &=
    \Sigma_{0,y}\otimes\left(P_x^{-1}\otimes I_y\right)\otimes E_{0y}
    \left(V-g\right) 
    + \Sigma_{Ny}\otimes\left(P_x^{-1}\otimes I_y\right)\otimes E_{N,y}
    \left(V-g\right )
    \\
    &\quad +\Sigma_{x0}\otimes\left(I_x\otimes P_y^{-1}\right)\otimes
    E_{x0} 
    \left(V-g\right)
    +\Sigma_{x,N}\otimes\left(I_x\otimes P_y^{-1}\right)\otimes E_{xN}
    \left(V-g\right),
  \end{aligned}
\end{equation}
\begin{theorem}
  \label{theo:dis}
  Assume that the velocity field $\U$ is a constant given by $\U =
  (u^1, u^2)^T$, and let $V$ be the semi-discrete solution defined by the
  scheme (\ref{eq:discrete}). Let $ u^{\ell,+} = \mx{u^{\ell}}{0} $ and $
  u^{\ell,-} = \mi{u^{\ell}}{0} $, for $\ell=1,\>2$. If  the penalty
  parameters satisfy
  \begin{equation}
    \label{eq:stacond}
    \begin{gathered}
      \sigma_1 \le -\frac {u^{1,+}}{2},\  
      \sigma_2 \le -\frac{u^{1,-}}{2},\
      \sigma_3 \le \frac {u^{2,+}}{2}\
      \quad {\rm and} \quad \sigma_4 \le \frac {u^{2,-}}{2}
    \end{gathered}
  \end{equation} 
  there exists positive constants $\alpha$ and $K$ such that
  \begin{equation}
    \label{eq:enestdisc}
    \norm{V(t)}^{2}\le \norm{\B_0}^2 + 
     \int\limits_{0}^t\int_{\partial\Omega} \g(t,x)\,dxd\tau,
 \end{equation}
and the scheme (\ref{eq:discrete}) is stable.
\end{theorem}
\begin{proof}
  The proof is similar to the standard way of proving stability of
  SBP-SAT schemes (see \cite{SN1}) and follows the proof for obtaining
  energy stability of the continuous problem in theorem \ref{theo:cp}.
  We outline the proof for the sake of completeness. For simplicity, we consider the case of constant velocities by setting $C = 0$ in (\ref{eq:discrete}) We start by multiplying (\ref{eq:discrete}) with $V^T(I_2\otimes
  P_x\otimes P_y)$ to obtain,
  \begin{equation}
    \label{eq:first}
    \begin{aligned}
      V^T&\left(I_2\otimes P_x\otimes P_y\right)\partial_t V
      \\
      &= -  V^T\left(\Lambda_x\otimes Q_x\otimes P_y\right)V 
      + V^T\left(\Lambda_y\otimes P_x\otimes Q_y\right)V 
      \\
      &\quad   +
      V^T(I_2\otimes P_x\otimes P_y) \otimes 
      \left[
        \begin{aligned}
          &\left(\Sigma_{0y}\otimes P_x^{-1}\otimes I_y\right) E_{0,y}
          + \left(\Sigma_{N,y}\otimes P_x^{-1}\otimes I_y\right) E_{N,y}  \\
        &\; +(\Sigma_{x,0}\otimes I_x\otimes P_y^{-1}) E_{x,0}
        +(\Sigma_{x,N}\otimes I_x\otimes P_y^{-1})E_{x,M}
        \end{aligned}
      \right]V.
    \end{aligned}
  \end{equation}
Adding this to its transpose and using the definition of SBP operators, we obtain
\begin{align*}
  \frac{d}{dt}&\norm V^2 \\
  &=
  -V^T\left(\Lambda_1\otimes \mathcal{B}_x \otimes P_y\right)V +
  V^T\left(\Lambda_2\otimes P_x\otimes \mathcal{B}_y\right)V 
  \\ 
  &\quad +
  2V^T(I_2\otimes P_x\otimes P_y)\otimes 
  \left[
    \begin{aligned}
      &\left(\Sigma_{0,y}\otimes P_x^{-1}\otimes I_y\right) E_{0,y}
      + \left(\Sigma_{N,y}\otimes P_x^{-1}\otimes I_y\right) E_{N,y}  \\
      &\; +\left(\Sigma_{x,0}\otimes I_x\otimes P_y^{-1}\right)
      E_{x,0} +\left(\Sigma_{x,N}\otimes I_x\otimes
        P_y^{-1}\right)E_{x,N}
    \end{aligned}
  \right]V,
\end{align*}
which implies
\begin{align*}
  \frac{d}{dt}\norm V^2 &=
  u^1\left(V^1_{0,y}\right)^T P_y \left(V^1_{0,y}\right) 
  - u^1\left(V^1_{N,y}\right)^T P_y \left(V^1_{N,y}\right)+
  u^1\left(V^2_{0,y}\right)^T P_y \left(V^2_{0,y}\right)
  \\ 
  &\quad
  - u^1\left(V^2_{N,y}\right)^T P_y \left(V^2_{N,y}\right)+
  u^2\left(V^1_{x,0}\right)^T P_x \left(V^1_{x,0}\right) 
  - u^2\left(V^1_{x,M}\right)^T P_x \left(V^1_{x,M}\right)
  \\
  &\quad 
  +u^2\left(V^2_{x,0}\right)^T P_x \left(V^2_{x,0}\right) 
  -u^2\left(V^2_{x,N}\right)^T P_x \left(V^2_{x,N}\right)
  \\
  &\quad + 2 \Bigl[
  \sigma_1 \left(\left(V^1_{0,y}\right)^T P_y \left(V^1_{0y}\right) 
  + \left(V^2_{0,y}\right)^T P_y \left(V^2_{0,y}\right)\right) 
  +\sigma_2 \left(V^1_{N,y}\right)^T P_y \left(V^1_{N,y}\right)
  \\
  &\hphantom{\quad + 2 \Bigl[}\quad
  + \sigma_2\left(V^2_{N,y}\right)^T P_y \left(V^2_{N,y}\right) 
  +\sigma_3 \left(\left(V^1_{x,0}\right)^T P_x \left(V^1_{x,0}\right)
  +\left(V^2_{x,0}\right)^T P_x \left(V^2_{x,0}\right)\right) 
  \\
  &\hphantom{\quad + 2 \Bigl[}\quad
  +\sigma_4 \left(\left(V^1_{x,N}\right)^T P_x \left(V^1_{x,N}\right)
  \left(V^2_{x,N}\right)^T P_x \left(V^2_{x,N}\right)\right)
  \Bigr].
\end{align*}
Using (\ref{eq:stacond}) and integrating in
time gives the energy estimate (\ref{eq:enestdisc}).
\end{proof}
\begin{remark}
  The above proof of stability assumes a constant velocity field. A
  proof of stability with a general velocity fields has been obtained
  in a recent paper \cite{svardsid3} by using the principle of frozen
  coefficients. The resulting stability estimate will lead to an exponential growth of energy (similar to (\ref{eq:enest})) due to the presence of lower order terms.
\end{remark}
We conclude this section with a few comments. For simplicity, we have only considered Cartesian meshes. However, the proofs are readily generalized to curvilinear grids by transforming the domain to a Cartesian. A stability proof is then obtained by freezing the coefficients. However, that requires $P$ to be diagonal, \cite{Svard04}. Furthermore, multi-block grids can also be handled and stable interfaces derived in a similar way as in, \cite{NordstromCarpenter01}.

\section{Numerical Experiments}
\label{sec:numex}

We test the SBP-SAT schemes of the previous section on a suite of
numerical experiments in order to demonstrate the effectiveness of
these schemes. We will use two different schemes : $SBP2$ and $SBP4$ scheme which are second-order (first-order) and fourth order (second-order) accurate in the interior (boundary) resulting in an overall second and third-order of accuracy. Time integration is performed by using a second order
accurate Runge-Kutta scheme at a $CFL$ number of $0.45$ for all numerical experiments. We found that using a fourth order
accurate Runge-Kutta scheme resulted in negligible differences in
the numerical results. The schemes have bounded errors, a typical behavior for hyperbolic equations with characteristic boundary conditions as shown in \cite{Nordstrom07}. Errors are propagated through the domain and leave the domain on account of the transparent boundaries. Hence, errors do not accumulate in time. On small domains, spatial errors become dominant. 

\vspace{0.5cm}
\noindent {\bf Numerical experiment $1$:}
In this experiment, we consider (\ref{eq:main}) with the
divergence-free velocity field $\U(x,y) = (-y,x)^T$. The
exact solution can be easily calculated by the method of
characteristics and takes the form
\begin{equation}
  \label{eq:ex}
  \B(\X,t)=R(t)\B_0(R(-t)\X),
\end{equation}
where $R(t)$ is a rotation matrix with angle $t$ and represents rotation of the initial data about the origin.

We consider the same test setup as in \cite{TF1} and \cite{fkrsid1} by
choosing the divergence free initial data,
\begin{equation}
  \B_0(x,y)=4
  \begin{pmatrix}
    -y\\ x-\frac{1}{2}
  \end{pmatrix}
  e^{-20\left((x-1/2)^2+y^2\right)},
\end{equation}
and the computational domain $[-1,1] \times [-1,1]$. Since the exact solution is known in
this case, one can in principle use this to specify the
boundary data $g$. Instead, we decided to
mimic a free space boundary (artificial boundary) by taking $g=0$. (which is a good guess at a far-field boundary).

\begin{figure}[htbp]
  \centering
\subfigure[half rotation, SBP2]{    \includegraphics[width=0.45\linewidth]{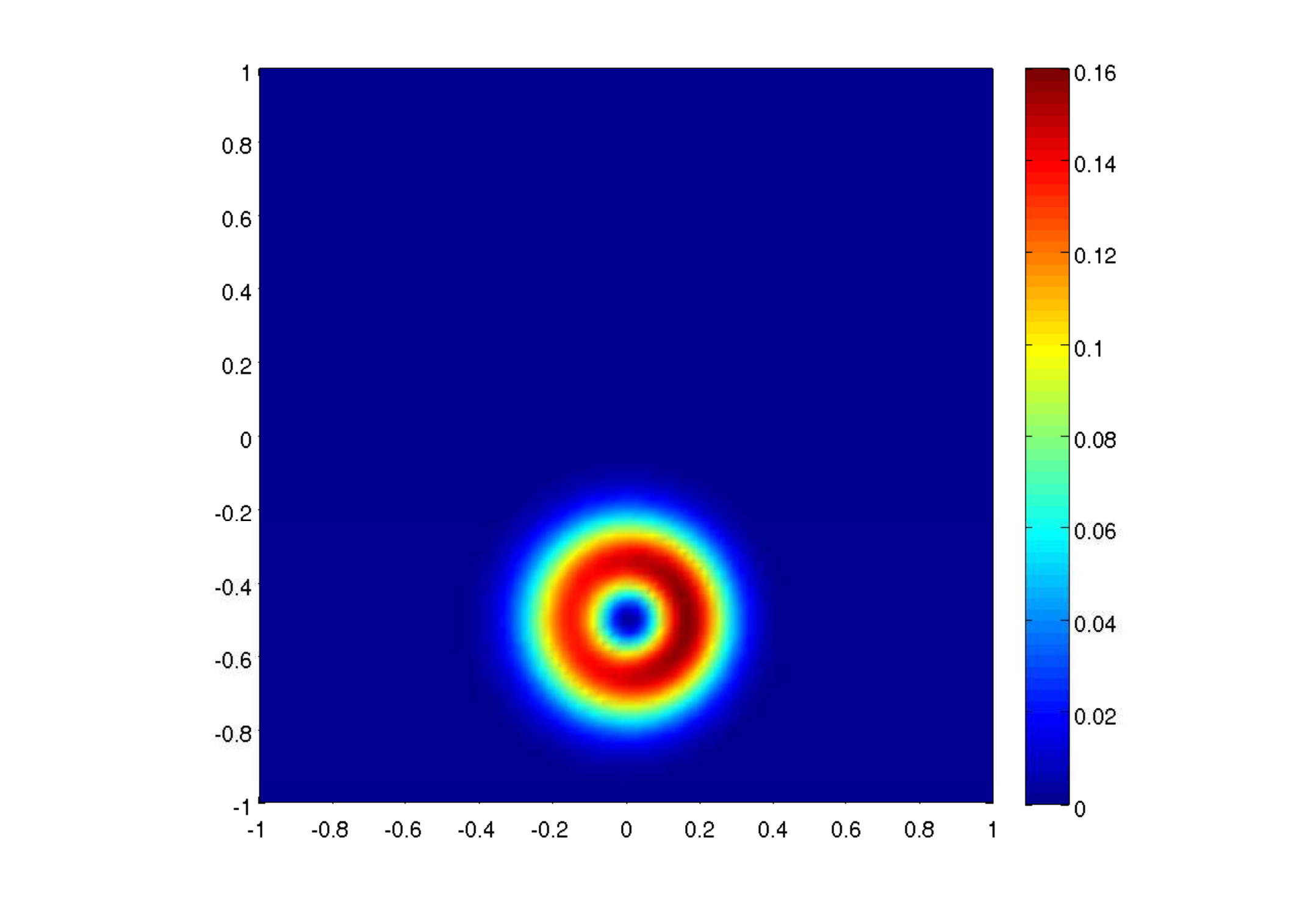} }
\subfigure[full rotation, SBP2]{    \includegraphics[width=0.45\linewidth]{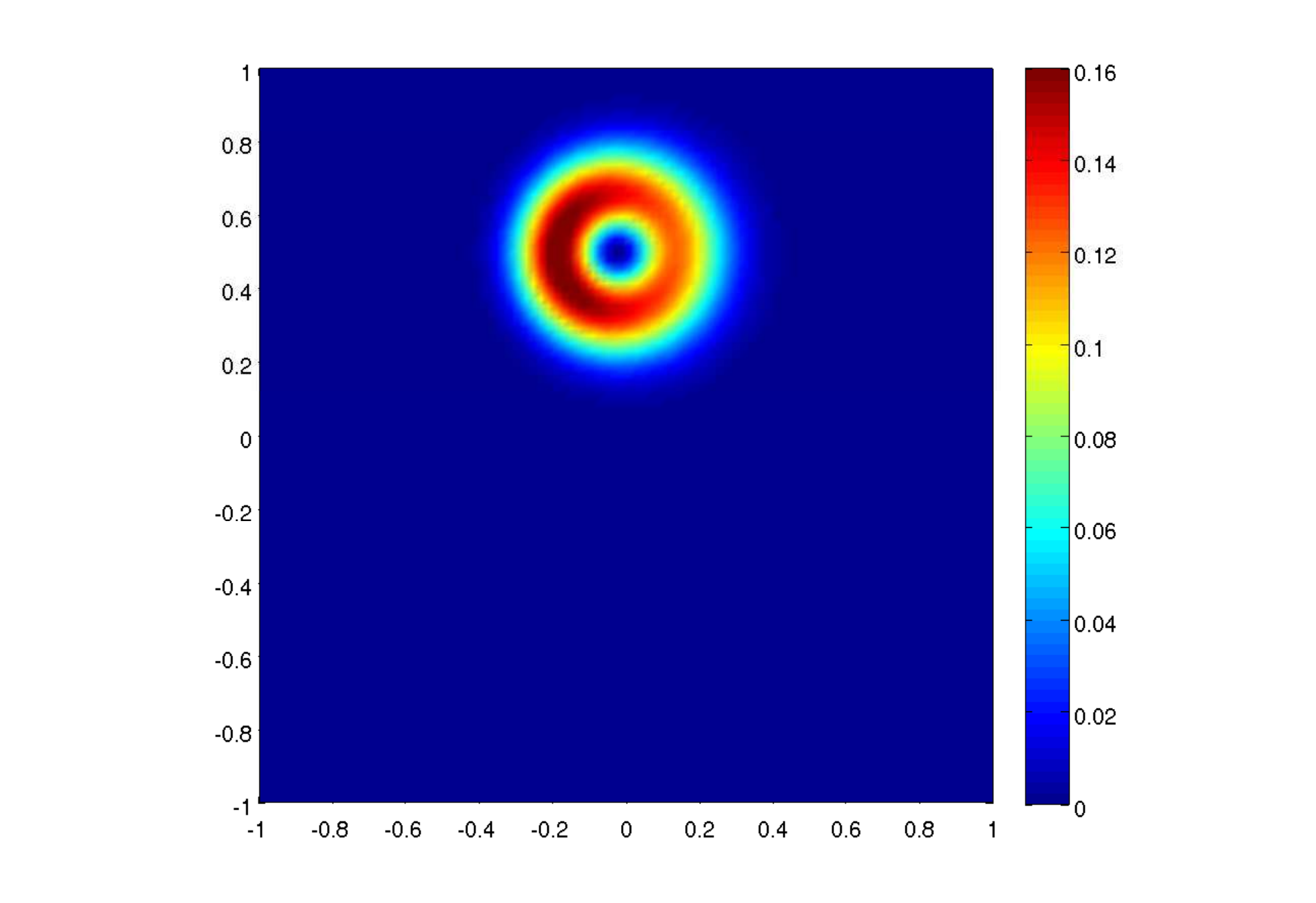}}
\subfigure[half rotation, SBP4]{    \includegraphics[width=0.45\linewidth]{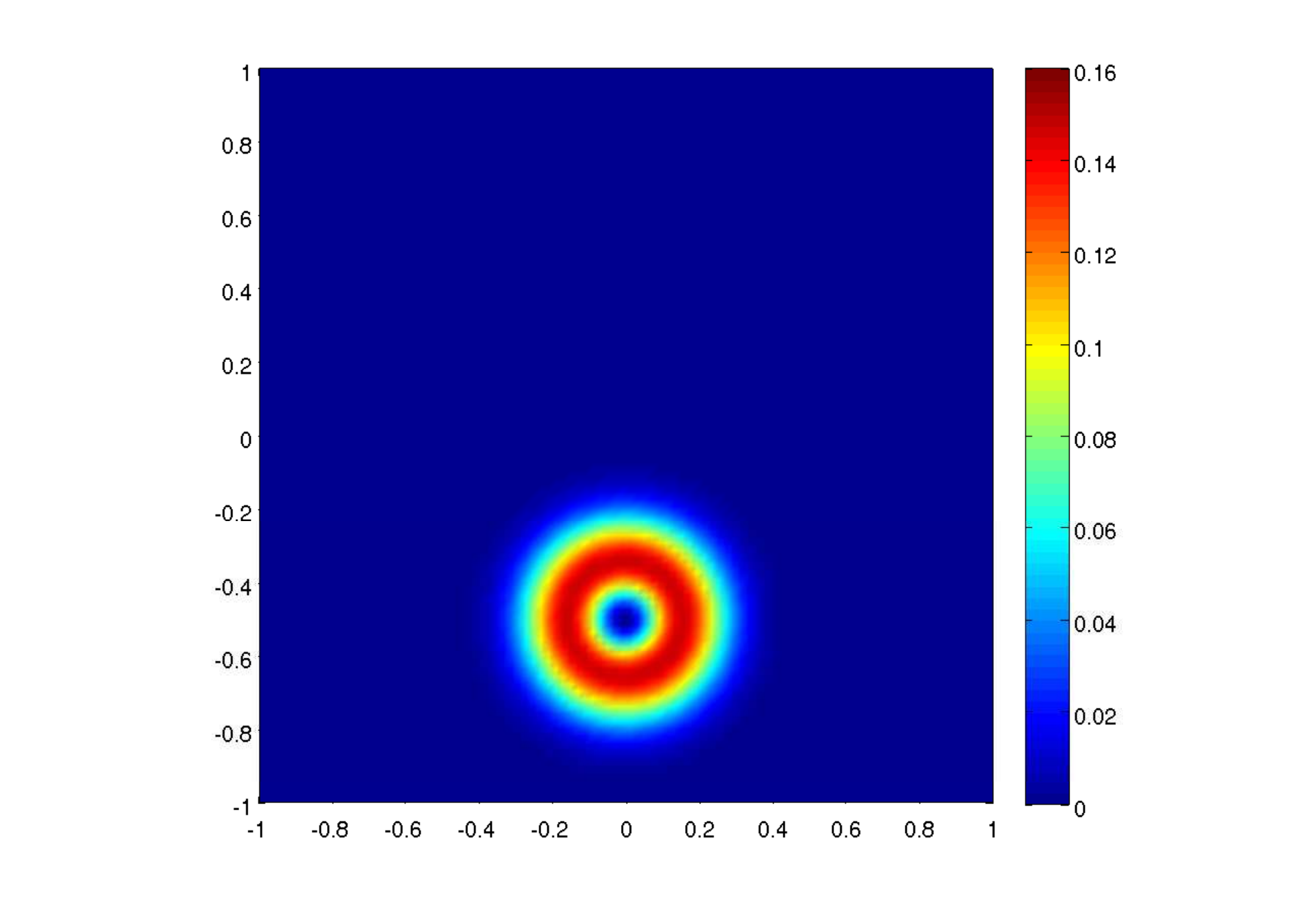}}
\subfigure[full rotation, SBP4]{    \includegraphics[width=0.45\linewidth]{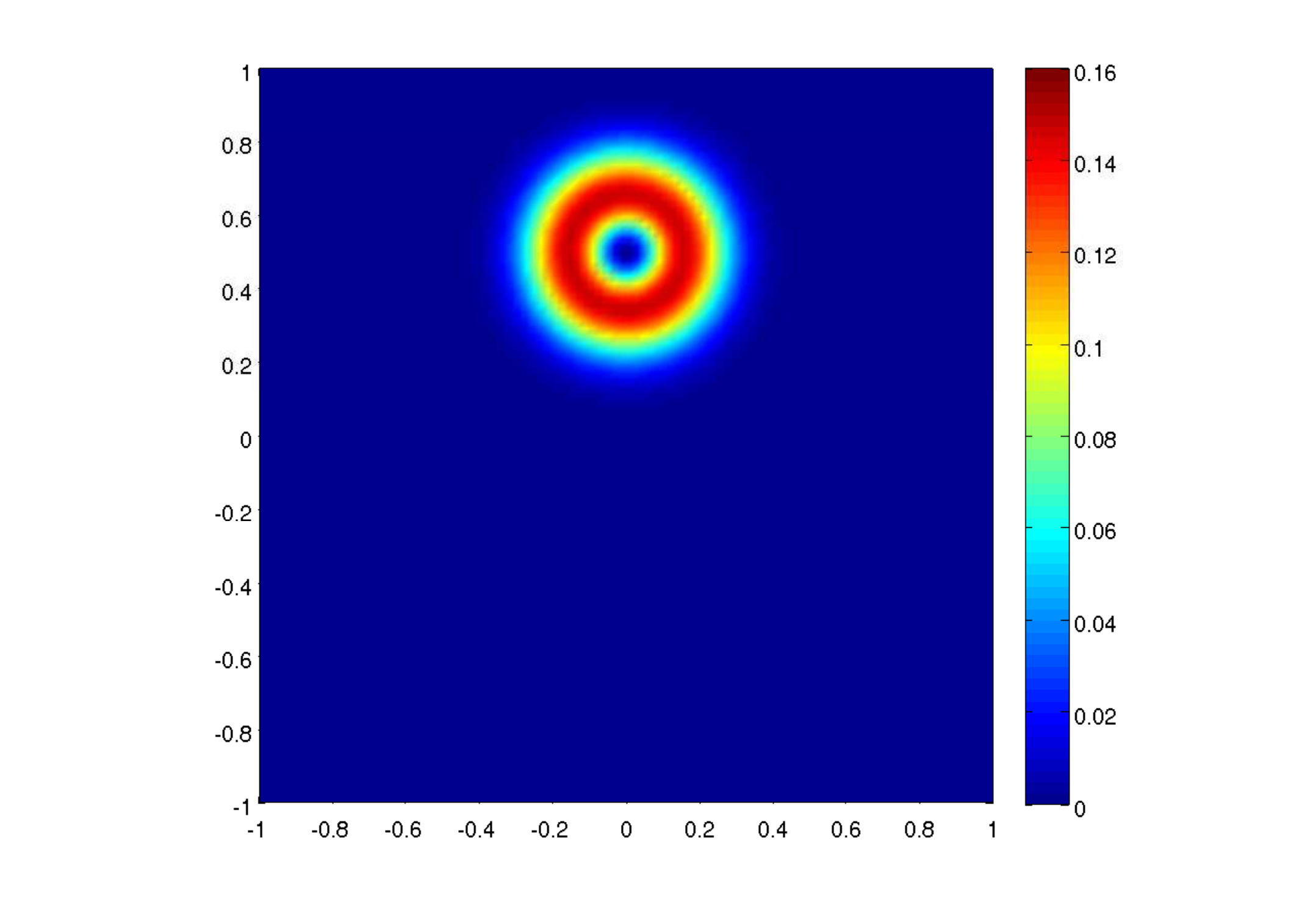} }
 
\caption{Numerical results for $|\B|$ in experiment $1$.}
\label{fig:1}  
\end{figure}

We run this test case with $SBP2$ and $SBP4$ schemes and present
different sets of results. In Figure \ref{fig:1}, we plot $|\B|=(|B^1|^2+|B^2|^2)^{1/2}$ at times $t= \pi$ (half-rotation) and $t=2\pi$ (one
full rotation) with the $SBP2$ and $SBP4$ schemes. 
As shown in this figure, $SBP2$ and $SBP4$ schemes resolve the
solution quite well. In fact, $SBP4$ is very accurate and keeps the
hump intact throughout the rotation.
\begin{table}[htbp]
\centering
\begin{tabular}{c|D{x}{}{-5} r D{x}{}{-5} r}
Grid size& \multicolumn{1}{c}{$SBP2$} 
& \multicolumn{1}{c}{rate} & 
\multicolumn{1}{c}{$SBP4$} & 
\multicolumn{1}{c}{rate}\\
\hline
40$\times$40   & 6.9x\e{1}     &      & 8.0x\e0   &    \\
80$\times$80   & 2.1x\e{1}     & 1.7  & 5.0x\e{-1}  & 4.0 \\
160$\times$160 & 5.5x\e{0}     & 2.0  & 4.5x\e{-2}  & 3.5 \\
320$\times$320 & 1.3x\e{0}     & 2.0  & 5.1x\e{-3}  & 3.1 \\
640$\times$640 & 3.3x\e{-1}    & 2.0  & 6.4x\e{-4}  & 3.0
\end{tabular}
\caption{Relative percentage errors in $l^2$ for $|\B|$ at time $t = 2
  \pi$ and rates of convergence for
  numerical experiment $1$ with $SBP2$ and $SBP4$ schemes.}
\label{tab:1}
\end{table}
In Table \ref{tab:1}, we present percentage relative errors in $l^2$. The
errors are computed at time $t = 2 \pi$ (one rotation) on a sequence
of meshes for both the $SBP2$ and $SBP4$ schemes. The results show
that the errors are quite low, particularly for $SBP4$ and the rate of convergence approaches the expected values of $2$ for $SBP2$ and $3$ for $SBP4$. Furthermore, the order of accuracy is
unaffected at these resolutions by using zero Dirichlet boundary data instead of the exact
solution at the boundary.

In order to compare the SBP schemes of this paper with other existing
schemes, we choose to compute the solutions for this problem with both the
first- and the second-order divergence-preserving scheme of~\cite{TF1}, which we label as
the $TF$ and $TF2$ schemes. Furthermore, we compute the solutions using the
first order stable upwind scheme designed in \cite{fkrsid1}, labeled 
 the $SUS$ scheme. The relative errors with each of these schemes are shown in
Table \ref{tab:6}.
\begin{table}[htbp]
\centering
\begin{tabular}{c|D{x}{}{-5} D{x}{}{-5} D{x}{}{-5} D{x}{}{-5} D{x}{}{-5}}
{Grid size}& 
\multicolumn{1}{c}{$SUS$} &
\multicolumn{1}{c}{$TF$} &
\multicolumn{1}{c}{$TF2$} 
\\
\hline
40$\times$40   &8.6x\e1 & 7.6x\e1  & 1.8x\e1         \\
80$\times$80   &7.3x\e1 & 6.4x\e1  & 1.3x\e1     \\
160$\times$160 &5.4x\e1 & 4.7x\e1  & 3.0x\e0      \\
320$\times$320 &3.6x\e1 & 3.3x\e1  & 1.0x\e0      \\
640$\times$640 &2.0x\e1 & 1.4x\e1  & 2.7x\e{-1}   
\end{tabular}
\caption{Relative percentage errors in $l^2$ for $|\B|$ at $t = 2
  \pi$ and for  numerical experiment $1$ with the $SUS$, $TF$, $TF2$, $SBP2$
  and the $SBP4$ schemes.}  
\label{tab:6}
\end{table} 
Results in Tables \ref{tab:1} and \ref{tab:6} show that the $TF$ and $SUS$ schemes lead
to similar errors and these errors are considerably larger than the
errors generated by the $TF2$ and $SBP2$ schemes, while the errors
generated by the $SBP4$ scheme are much smaller again.  

A fair comparison of the the five schemes SUS, TF, TF2, SBP2 and SBP4 requires information on the computational work with each scheme for the same error level. We observe from tables \ref{tab:1} and \ref{tab:6} that for a given relative error of approximately $20$ percent, the first-order SUS scheme requires a $640 \times 640$ mesh, the TF scheme requires a $500 \times 500$ mesh (based on extrapolation from table \ref{tab:6}), whereas both the second-order schemes require meshes coarser than a $50 \times 50$ mesh. The fourth-order scheme yields similar error levels on even coarser meshes. Thus, the second-order schemes require about $1\%$ of the total grid points to the first-order schemes to produce comparable errors.  Even taking into account that the second order schemes use more  
operations for each grid point, this still makes the
second order schemes at least $25-30$ times more efficient than the first order  schemes. Similarly an error level of about one percent is attained with SBP2 on a $320 \times 320$ mesh, with TF2 on a similar $320 \times 320$ mesh and with SBP4 on a $50 \times 50$ mesh.  Thus the second order schemes need about $36$ times more grid points to produce errors similar to those of the fourth order schemes. Taking extra work for the fourth-order scheme per grid point into account, we still get that the fourth-order scheme is roughly $10$ times more efficient than the second-order schemes. These numbers are approximations but display a clear qualitative trend i.e., it is much more efficient to use high-order schemes for the induction equations.

As the solution (\ref{eq:ex}) in this case is smooth, it is also a
solution for the constrained form (\ref{eq:Maxwell3D}). Furthermore,
the initial data is divergence free and so is the exact solution. We
did not attempt to preserve any particular form of discrete divergence
while designing the SBP schemes (\ref{eq:discrete}). A natural thing would be show that some form of discrete divergence produced by the schemes was
bounded in $l^2$. We were unable to obtain such a divergence bound for
(\ref{eq:discrete}) in this paper. A related SBP-SAT scheme for the
``conservative'' symmetric form (\ref{eq:induc2}) with SBP operators
for discretizing spatial derivatives coupled with a novel
discretization of the source terms in (\ref{eq:induc2}) was shown to
have bounded discrete divergence in a recent paper \cite{svardsid3}.

In the absence of a rigorous divergence bound, we proceed to examine
how divergence errors generated by the SBP schemes behave and whether
they have any impact on the quality of the discretization. We define
the following discrete divergence,
\begin{align*}
  \mathrm{div}_P(V)= (P_x^{-1}Q_x\otimes I_{y})V^{1} + (I_{x}\otimes
  P_y^{-1}Q_y)V^{2}.
\end{align*}
This corresponds to the standard centered discrete divergence operator
at the corresponding order of accuracy.  The divergence errors in
$l^2$ and rates of convergence at time $t=2\pi$ for the $SBP2$ and
$SBP4$ schemes on a sequence of meshes are presented in
Table~\ref{tab:3}.
\begin{table}[htbp]
  \centering
  \begin{tabular}{c|D{x}{}{-5} r D{x}{}{-5} r}
    Grid size& \multicolumn{1}{c}{$SBP2$} 
    & \multicolumn{1}{c}{rate} & 
    \multicolumn{1}{c}{$SBP4$} & 
    \multicolumn{1}{c}{rate}\\
    \hline
    20$\times$20  & 1.0x\e0    &     & 7.3x\e{-1}   &     \\
    40$\times$40  & 8.0x\e{-1} & 0.4 & 1.2x\e{-1}   & 2.6 \\
    80$\times$80  & 2.7x\e{-1} & 1.6 & 8.2x\e{-3}   & 3.8 \\
    160$\times$160& 7.0x\e{-2} & 2.0 & 1.0x\e{-3}   & 3.0 \\
    320$\times$320& 2.5x\e{-2} & 1.5 & 1.7x\e{-4}   & 2.6
  \end{tabular}
  \caption{Numerical Experiment $1$:
    Divergence (errors) in $l^2$ and rates of convergence at time
    $t=2\pi$ for both the $SBP2$ and $SBP4$ schemes. }
  \label{tab:3} 
\end{table}
From Table~\ref{tab:3}, we conclude that although the initial
divergence is zero, the discrete divergence computed with both the
$SBP2$ and $SBP4$ schemes is not zero. However, the divergence errors
are very small even on fairly coarse meshes and converge to zero at a
rate of $1.5$ and $2.5$ for $SBP2$ and $SBP4$ schemes respectively. A
simple truncation error analysis suggests that these rates for the $SBP2$
and $SBP4$ schemes are optimal. The quality of the approximations is 
good and the rates of convergence do not seem to suffer from  not preserving
any form of discrete divergence.

In order to compare with existing schemes, we compare the divergence
errors generated by the $SUS$, $TF$ and the $TF2$ schemes with the
$SBP2$ and the $SBP4$ schemes in table~\ref{tab:7}. 
\begin{table}[htbp]
  \centering
  \begin{tabular}{c | D{x}{}{5} D{x}{}{5} D{x}{}{5} D{x}{}{5} D{x}{}{5}}
    {Grid size}& 
    \multicolumn{1}{c}{$SUS$} &
    \multicolumn{1}{c}{$TF$} &
    \multicolumn{1}{c}{$TF2$} &
     \\
    \hline
    40$\times$40  &1.1x\e{-1}&2.7x\e{-2}&1.2x\e{-2}\\
    80$\times$80  &1.3x\e{-1}&1.7x\e{-2}&4.0x\e{-3}\\
    160$\times$160&1.4x\e{-1}&1.4x\e{-2}&2.4x\e{-3}\\
    320$\times$320&1.1x\e{-1}&1.2x\e{-2}&9.7x\e{-4}
  \end{tabular}
  \caption{Numerical Experiment $1$: The discrete divergence
    $\mathrm{div}_P$ in $l^2$ at 
    $t=2\pi$ for the $SUS$, $TF$ and $TF2$   schemes.}
  \label{tab:7} 
\end{table}
From Table~\ref{tab:7}, we can draw the following conclusions about
divergence errors. The $SUS$ scheme is not tailored to preserve any form
of discrete divergence. The divergence errors generated by this scheme
seems to be low on coarse meshes. The $TF$ and $TF2$ schemes are designed to preserve a special form of
discrete divergence which is different from the standard central form.
Nevertheless, the analysis presented in~\cite{TF1} suggested that the errors in
the standard divergence operator will also be quite low. This is
indeed the case. On the coarser meshes, the divergence is much larger for
the $SBP2$ scheme than the $TF$ schemes, but from Table~\ref{tab:6} we
see that the errors in the solution are similar.

Furthermore, the divergence errors converge quickly for the
 $SBP4$ scheme, as well as as the for the $TF2$
scheme. The above results indicate that controlling some form of discrete
divergence is not necessary to approximate solutions of the magnetic
induction equations in a stable and accurate manner.

\begin{figure}[htbp]
  \centering

\subfigure[$t=0$]{    \includegraphics[width=0.32\linewidth]{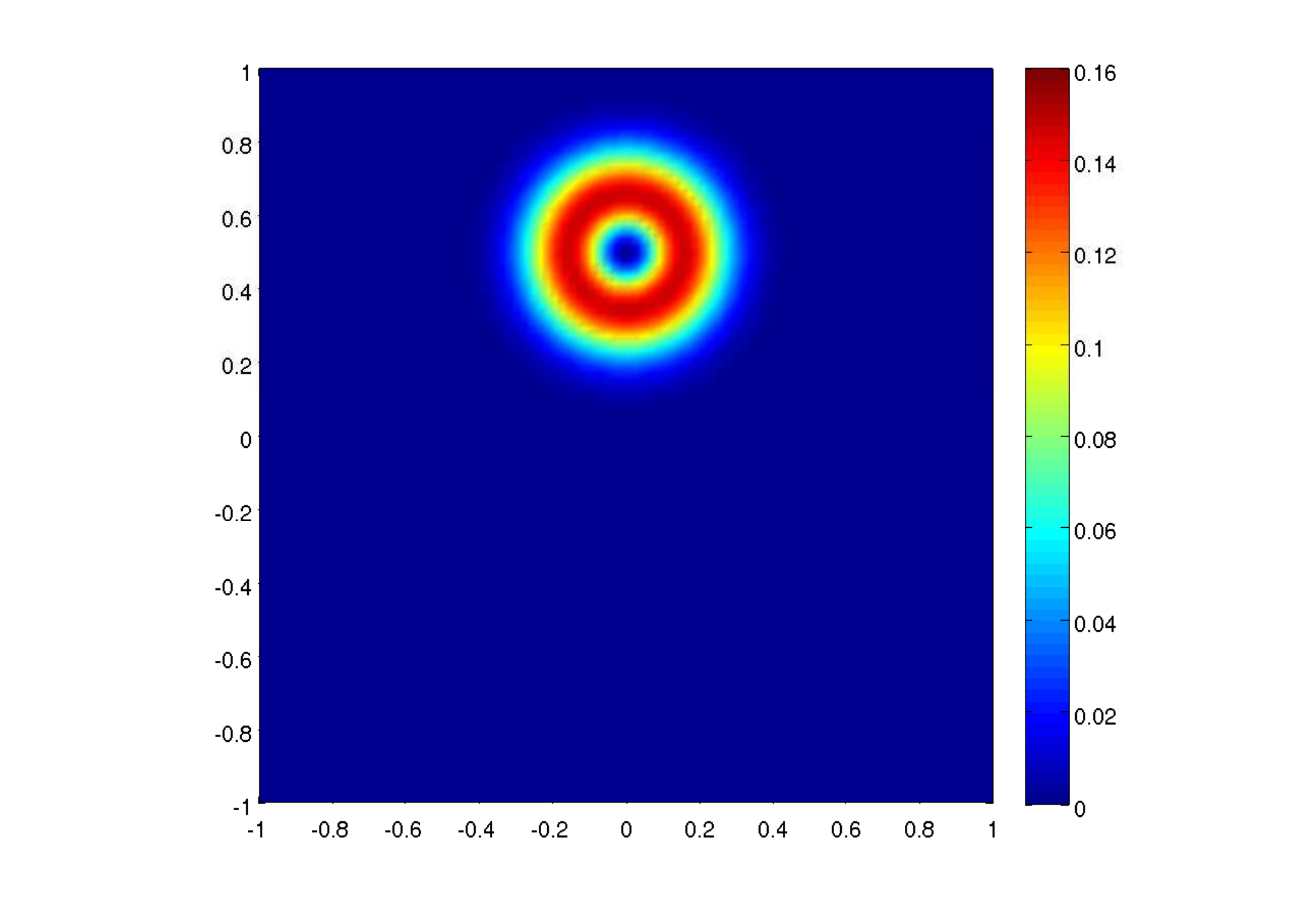} }
\subfigure[$t=10\pi$, SBP2]{    \includegraphics[width=0.32\linewidth]{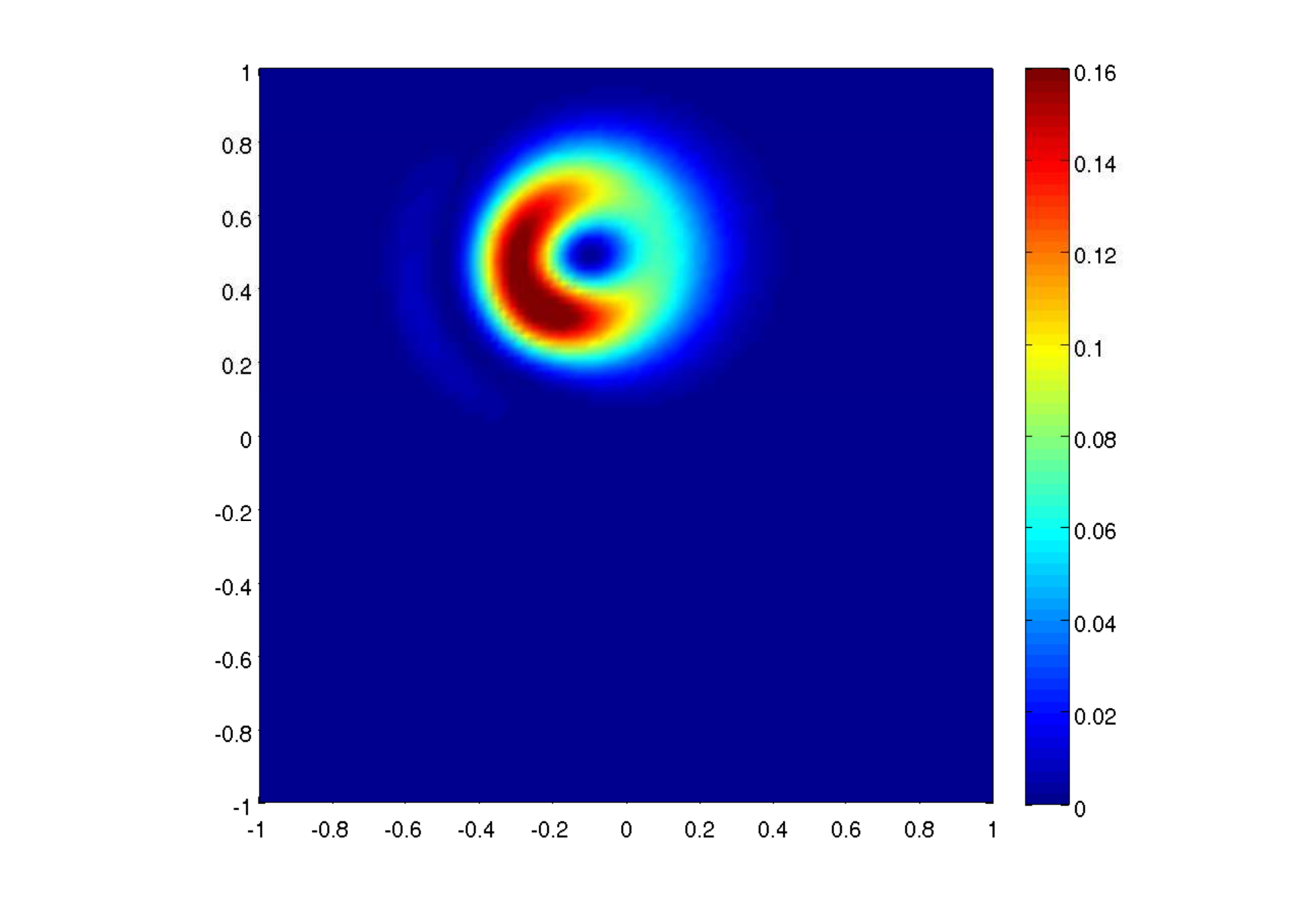}}
\subfigure[$t=100\pi$, SBP4]{    \includegraphics[width=0.32\linewidth]{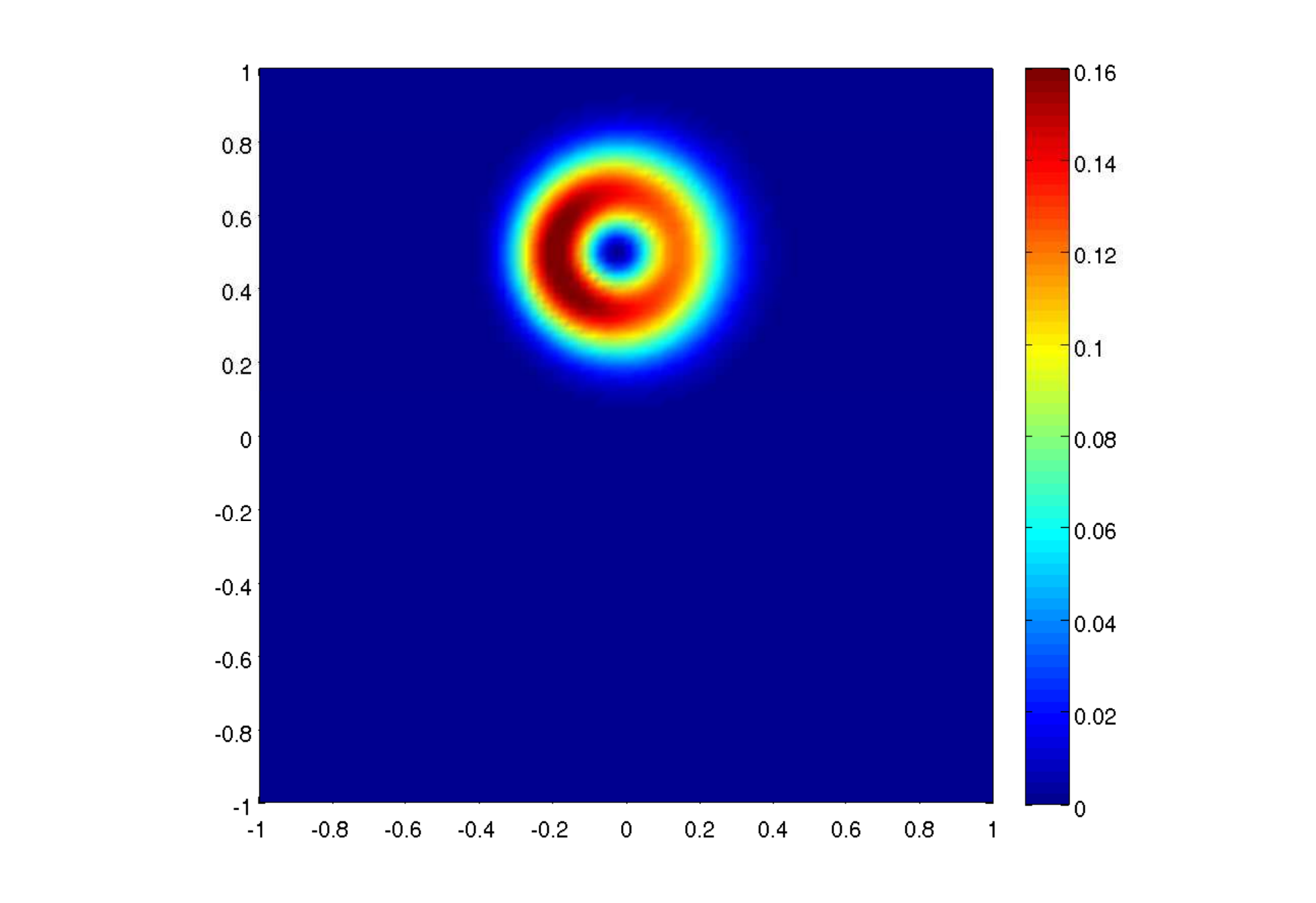}}
  \caption{ Numerical results for $|\B|$ in experiment $1$.}
  \label{fig:2}
\end{figure}

Next, we consider long time integration. The energy
estimate (\ref{eq:enestdisc}) suggests that the energy of the
approximate solutions can grow exponentially in time. In order to test this we computed approximate
solutions with the $SBP2$, $SBP4$  and the $TF2$ schemes till time $t = 100 \pi$, i.e., for fifty full rotations on a $100 \times 100$ mesh.  The numerical results in are presented
in Figure~\ref{fig:2} and Table~\ref{tab:2}.
\begin{table}[htbp]
  \centering
  \begin{tabular}{c | D{x}{}{-5} D{x}{}{-5} D{x}{}{-5}}
    $2\pi t$ & 
    \multicolumn{1}{c}{$SBP2$} &
    \multicolumn{1}{c}{$SBP4$} &
    \multicolumn{1}{c}{$TF2$} \\
    \hline
    $t=1$    & 2.1x\e{1} & 5.1x\e{-1} & 8.8x\e0\\
    $t=5$    & 7.7x\e{1} & 2.7x\e{0}  & 3.2x\e1\\
    $t=10$   & 1.0x\e{2} & 4.7x\e{0}  & 5.0x\e1\\
    $t=15$   & 1.1x\e{2} & 6.6x\e{0}  & 6.3x\e1\\
    $t=20$   & 1.2x\e{2} & 8.7x\e{0}  & 7.2x\e1\\
    $t=30$   & 1.2x\e{2} & 1.9x\e{1}  & 8.4x\e1\\
    $t=40$   & 1.3x\e{2} & 3.1x\e{1}  & 9.2x\e1\\
    $t=50$   & 1.4x\e{2} & 4.3x\e{1}  & 1.0x\e2
  \end{tabular}
  \caption{Relative percentage  $l^2$ errors in $| \B |$ with
    $SBP2$, $SBP4$ and $TF2$ for numerical experiment $1$.}
  \label{tab:2} 
\end{table}
These computations were performed on a fixed $100 \times 100$ mesh. In
Figure~\ref{fig:2}, we compare the $SBP2$ and $SBP4$ schemes after
five and fifty rotations respectively. We see that after 5 rotations, 
$SBP2$ gives a ``hump'' which is somewhat smeared and with a
pronounced asymmetry. 
 On the other hand, the hump produced by the 
$SBP4$ scheme is much more accurate.
As shown in Table \ref{tab:2}, the absolute errors with the $SBP4$ scheme are much
lower than the errors due to the second-order schemes $SBP2$ and
$TF2$. In fact, the errors with $SBP2$ after just five rotations are
about three times the error with $SBP4$ after fifty rotations. This
experiment makes a strong case for using high-order schemes for 
problems requiring long time integration.

\vspace{0.5cm}
\noindent {\bf Numerical Experiment $2$:}
In the previous numerical experiment, the hump was confined to the
interior of the domain during the rotation. Hence, the choice of zero
Dirichlet data at the boundary was reasonable and led to stable and
accurate approximations. In order to illustrate the effect of the
boundary better, we choose the computational domain $[0,1] \times [0,1]$ and use the same velocity field and initial data as in the previous experiment. Now, the hump ``exits'' the domain at one part of the boundary (including a corner)
and will re-enter the domain from another part of the boundary. The
choice of boundary discretization becomes crucial in this case. 

We select the exact solution (\ref{eq:ex}) restricted to the boundary
as the Dirichlet boundary data in (\ref{eq:discrete}).
\begin{figure}[htbp]
  \centering

   \subfigure[SBP2, $t=\pi/2$]{ \includegraphics[width=0.45\linewidth]{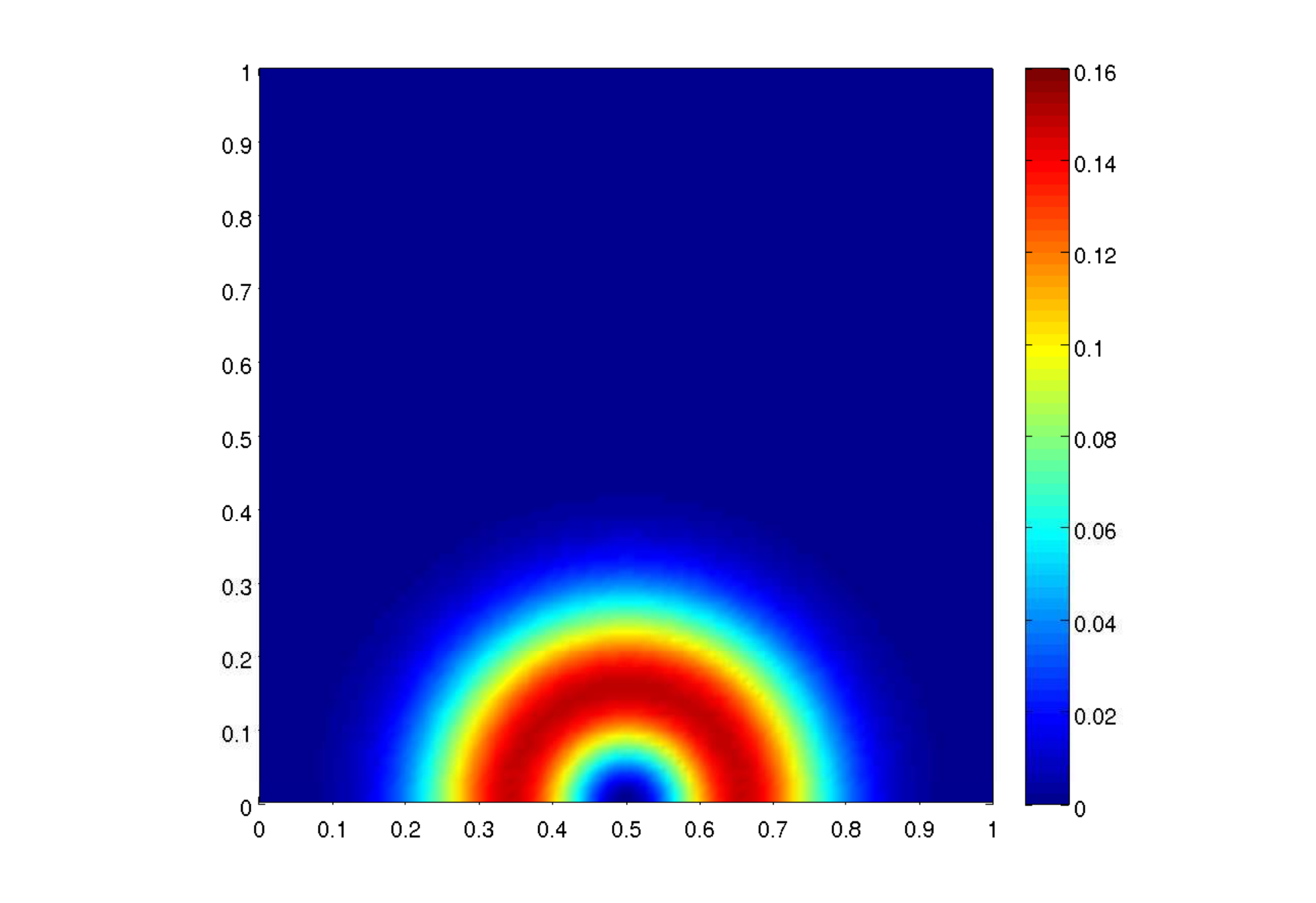} }
\subfigure[SBP2, $t=2\pi$]{    \includegraphics[width=0.45\linewidth]{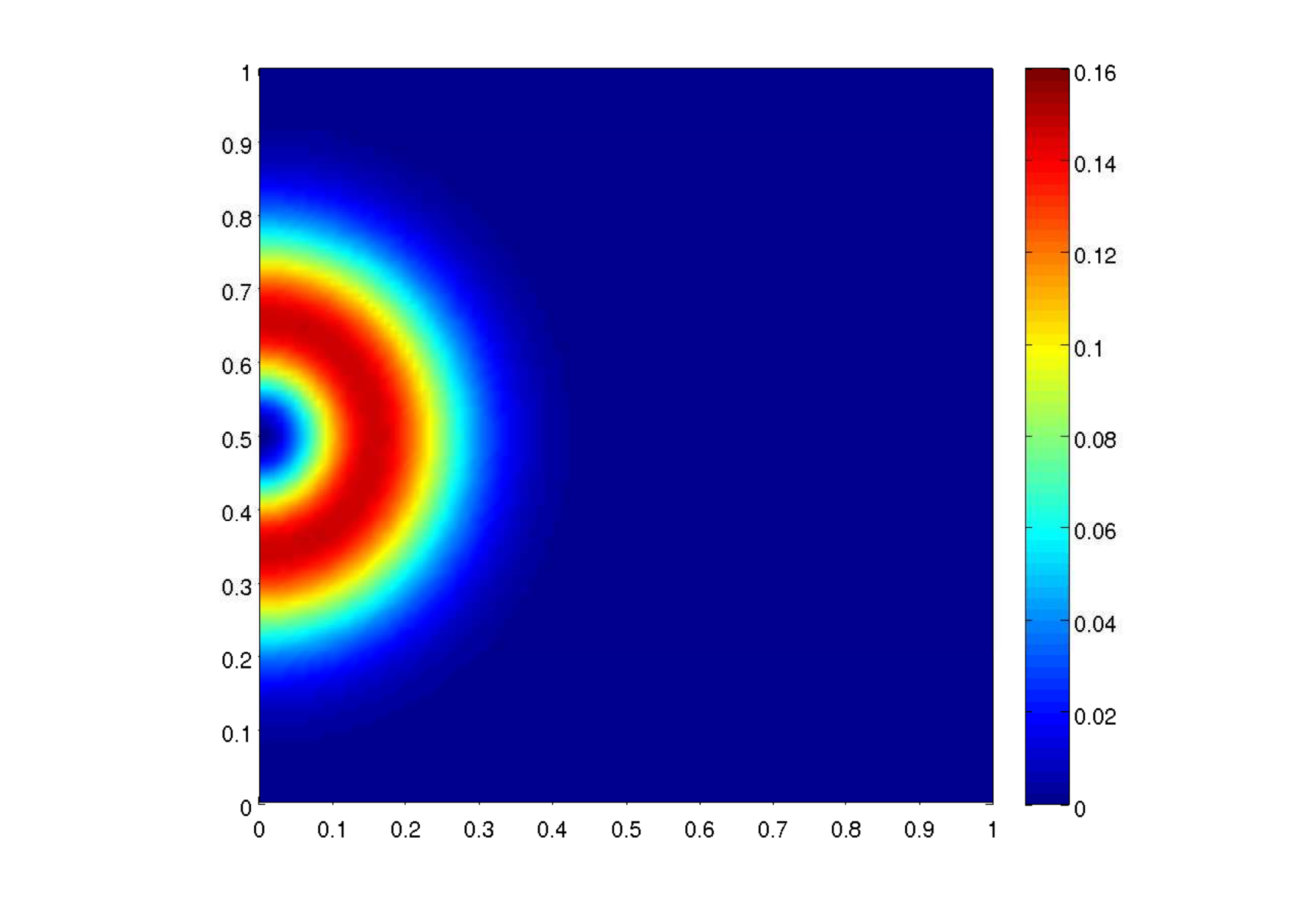}}

\subfigure[SBP4,$t=\pi/2$]{    \includegraphics[width=0.45\linewidth]{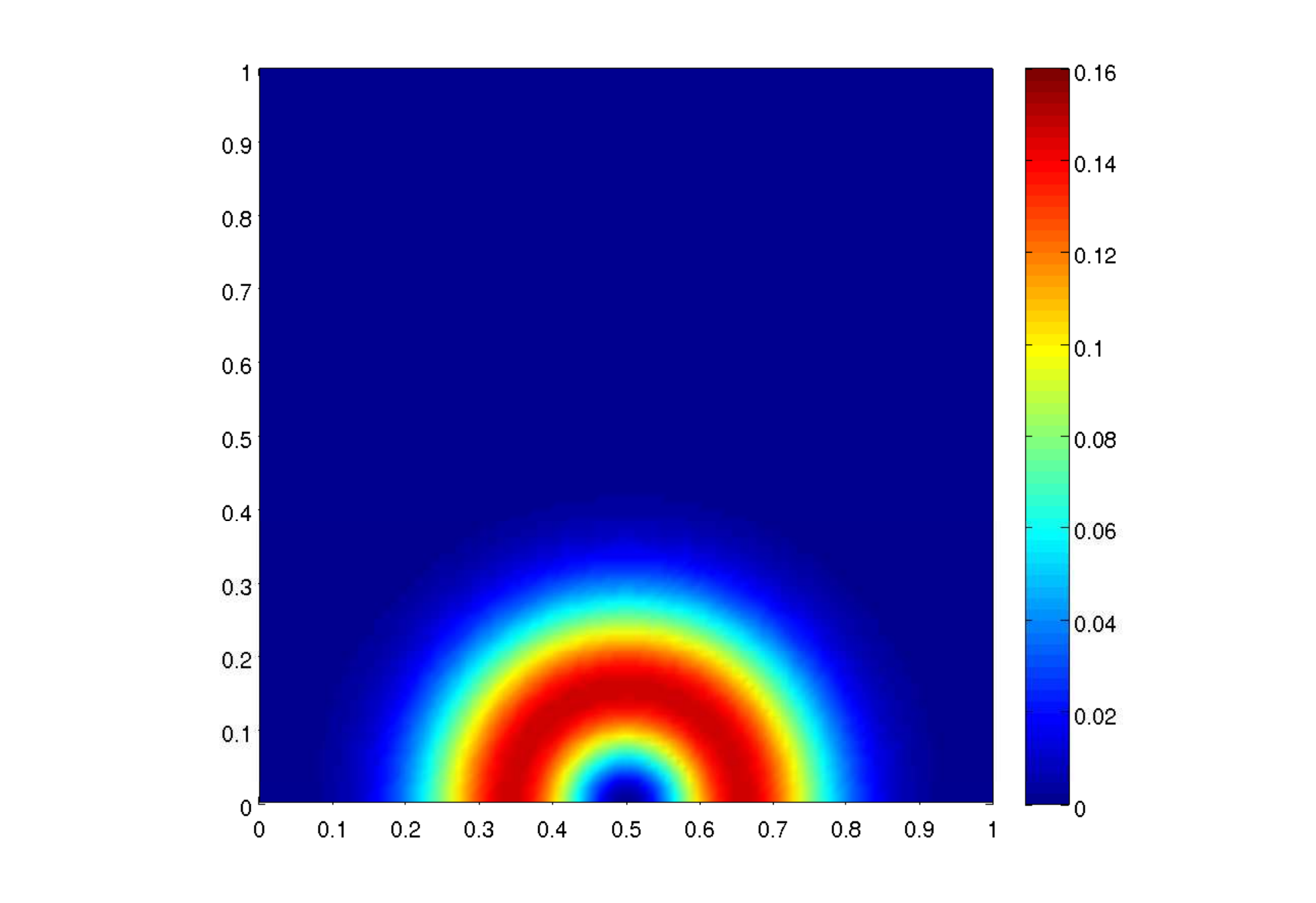}}
\subfigure[SBP4,$t=2\pi$]{    \includegraphics[width=0.45\linewidth]{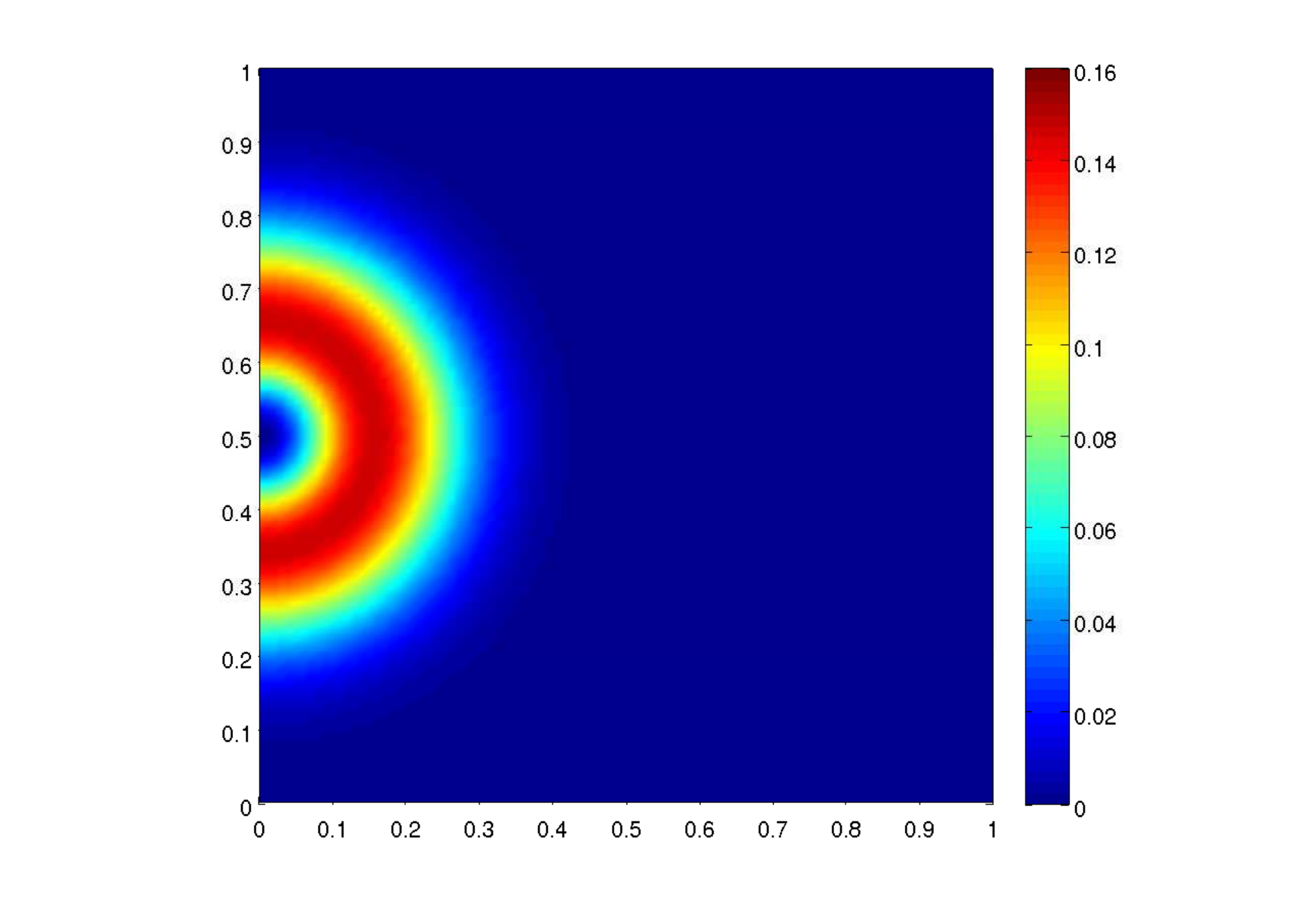}}

  \caption{Numerical results for experiment $2$. Mesh  size $100\times 100$.}

  \label{fig:4} 
\end{figure}
In Figure~\ref{fig:4}, the approximate solutions computed with both
$SBP2$ and $SBP4$ on a $100 \times 100$ mesh are plotted at time $t=
\pi/2$ (quarter rotation) and time $t = 2 \pi$ (full rotation). As shown in
this figure, both schemes perform very well. The hump at both the exit as well as the re-entry is 
clearly resolved with no noticeable numerical artefacts or
reflections. 
\begin{table}[htbp]
  \centering
  \begin{tabular}{c|D{x}{}{-5} r D{x}{}{-5} r}
    Grid size& \multicolumn{1}{c}{$SBP2$} 
    & \multicolumn{1}{c}{rate} & 
    \multicolumn{1}{c}{$SBP4$} & 
    \multicolumn{1}{c}{rate}\\
    \hline
    10$\times$10  & 2.5x\e{1}  &     & 1.1x\e{1}  &     \\
    20$\times$20  & 5.8x\e{0}  & 2.1 & 1.5x\e{0}  & 2.9 \\
    40$\times$40  & 1.3x\e{0}  & 2.0 & 1.6x\e{-1} & 3.3 \\
    80$\times$80  & 3.0x\e{-1} & 2.0 & 1.6x\e{-2} & 3.2 \\
    160$\times$160& 7.4x\e{-2} & 2.0 & 1.9x\e{-3} & 3.1  
  \end{tabular}
  \caption{Numerical experiment $2$: Relative percentage errors for $|\B|$ in $l^2$
    and rates of convergence for both 
    $SBP2$ and $SBP4$.} 
  \label{tab:4}
\end{table}

\begin{table}[htbp]
  \centering
  \begin{tabular}{c|D{x}{}{-5} r D{x}{}{-5} r}
    Grid size& \multicolumn{1}{c}{$SBP2$} 
    & \multicolumn{1}{c}{rate} & 
    \multicolumn{1}{c}{$SBP4$} & 
    \multicolumn{1}{c}{rate}\\
    \hline
    10$\times$10  & 6.4x\e{-1}  &     & 9.7x\e{-2}  &     \\
    20$\times$20  & 3.9x\e{-1}  & 0.7 & 2.4x\e{-2}  & 2.0 \\
    40$\times$40  & 9.1x\e{-2}  & 2.2 & 1.9x\e{-3} & 3.6 \\
    80$\times$80  & 2.6x\e{-2} & 1.8  & 3.0x\e{-4} & 2.7 \\
    160$\times$160& 8.9x\e{-3} & 1.6  & 5.1x\e{-5} & 2.5  
  \end{tabular}
  \caption{Numerical experiment $2$: Divergence (errors) in $l^2$
    and rates of convergence for both $SBP2$ and $SBP4$ at time $t=2\pi$.} 
  \label{tab:4b}
\end{table}

As shown in Table \ref{tab:4}, the errors are low after one full
rotation for both the $SBP2$ and $SBP4$ schemes. In fact, the size of
relative errors is lower than in the previous numerical experiment. As
expected, the rates of convergence tend to $2$ and $3$ for $SBP2$ and
$SBP4$ respectively. In Table \ref{tab:4b} the divergence errors and their convergence rates are listed. They are small and the convergences approach the expected values $1.5$ and $2.5$.

On the other hand, when we tried to compute this example with the divergence preserving $TF$ and $TF2$ schemes, the solution blew up on account of boundary instabilities.

\vspace{0.5cm}
\noindent {\bf Numerical Experiment 3: (Discontinuous solutions.)}
As remarked earlier, the
magnetic induction equations (\ref{eq:main}) are a sub-model in the
nonlinear MHD equations. As a consequence, one must solve the
induction equation with both
discontinuous velocity fields and initial data. It is
therefore interesting to see how well the SBP-SAT schemes handle
discontinuous velocity fields and initial data. 

The SBP operators use centered finite differences in the interior.
It is well known that using central differences leads to oscillations
around discontinuities. Therefore the SBP schemes cannot be used directly in
this regime, see~\cite{MSN1} for details. To calculate solutions with discontinuities, one adds a small
amount of 
explicit numerical diffusion that retain the accuracy of the first
derivative SBP approximations as well as maintain the energy stability
of the SBP scheme. We will use these
operators together with the $SBP2$ and $SBP4$ schemes in order to
compute discontinuous solutions of the magnetic induction equations.

The second-order (fourth-order) SBP operator for the first derivative with a second-order (fourth-order) numerical diffusion operator gives an
approximation which is formally second-order (fourth-order) accurate in the interior
of the computational domain. It turns out that a different scaling
(dividing by the mesh size) of the numerical diffusion operator leads
to a first order (third-order) ``upwind'' scheme. We will test all these numerical
diffusion operators a  numerical experiment first described in
\cite{fkrsid1}.

The computational domain is $[0,1] \times [0,1]$. Consider the
constant velocity field, $\U = (1,2)^T$ and the discontinuous initial
data,
$$
B^1_0(x,y)=B^2_0(x,y)=
\begin{cases}
  2 &\text{if $x>y$},\\
  0 &\text{otherwise.}
\end{cases}
$$
In this case, the exact solution (see \cite{fkrsid1}) of
(\ref{eq:main}) reads
$$
\B(x,y,t)=\B_0(x-t,y-2t).
$$
The initial discontinuity simply moves along the diagonal of
the domain.  We use the exact solution restricted to the boundary as
the Dirichlet boundary data. Tests with generic $SBP$-$SAT$
schemes, (\ref{eq:discrete}), showed that the approximate solutions were
very oscillatory, and we damp these oscillations by adding numerical
diffusion.

\begin{figure}[htbp]
  \centering

\subfigure[$SBP1$]{\includegraphics[width=0.45\linewidth]{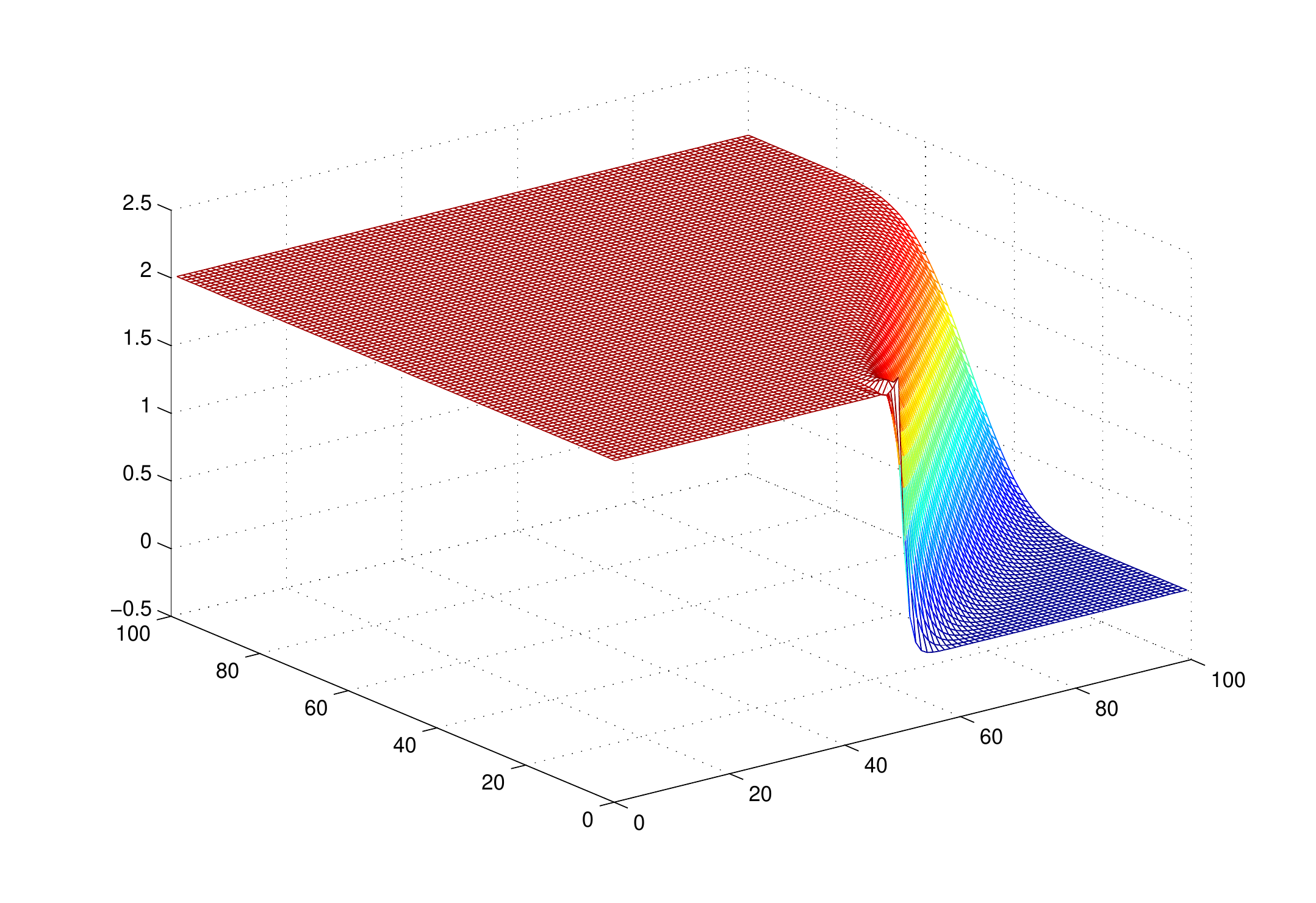} }
\subfigure[$SBP2$, second-order diffusion]{\includegraphics[width=0.45\linewidth]{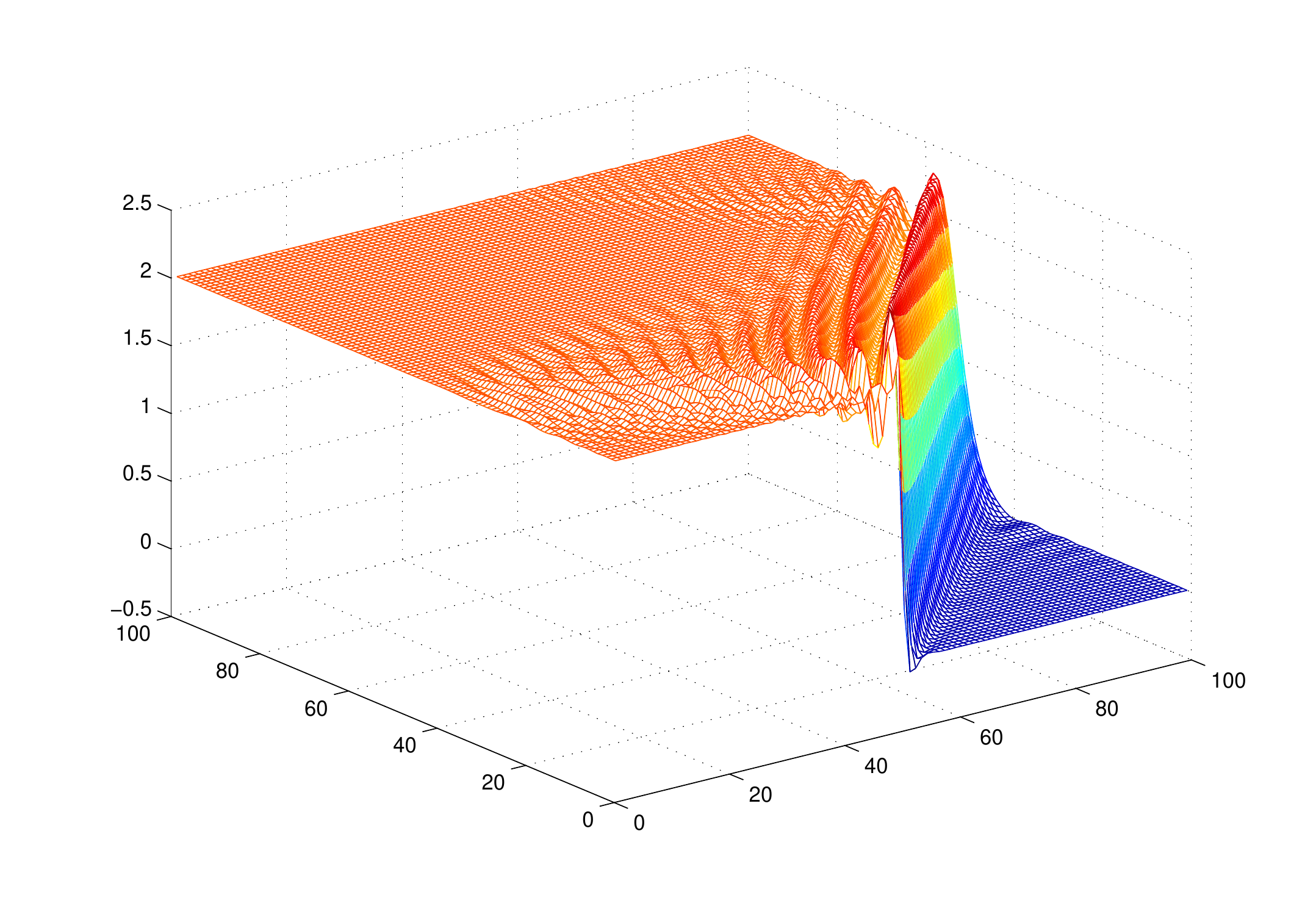} }
\subfigure[$SBP3$]{\includegraphics[width=0.45\linewidth]{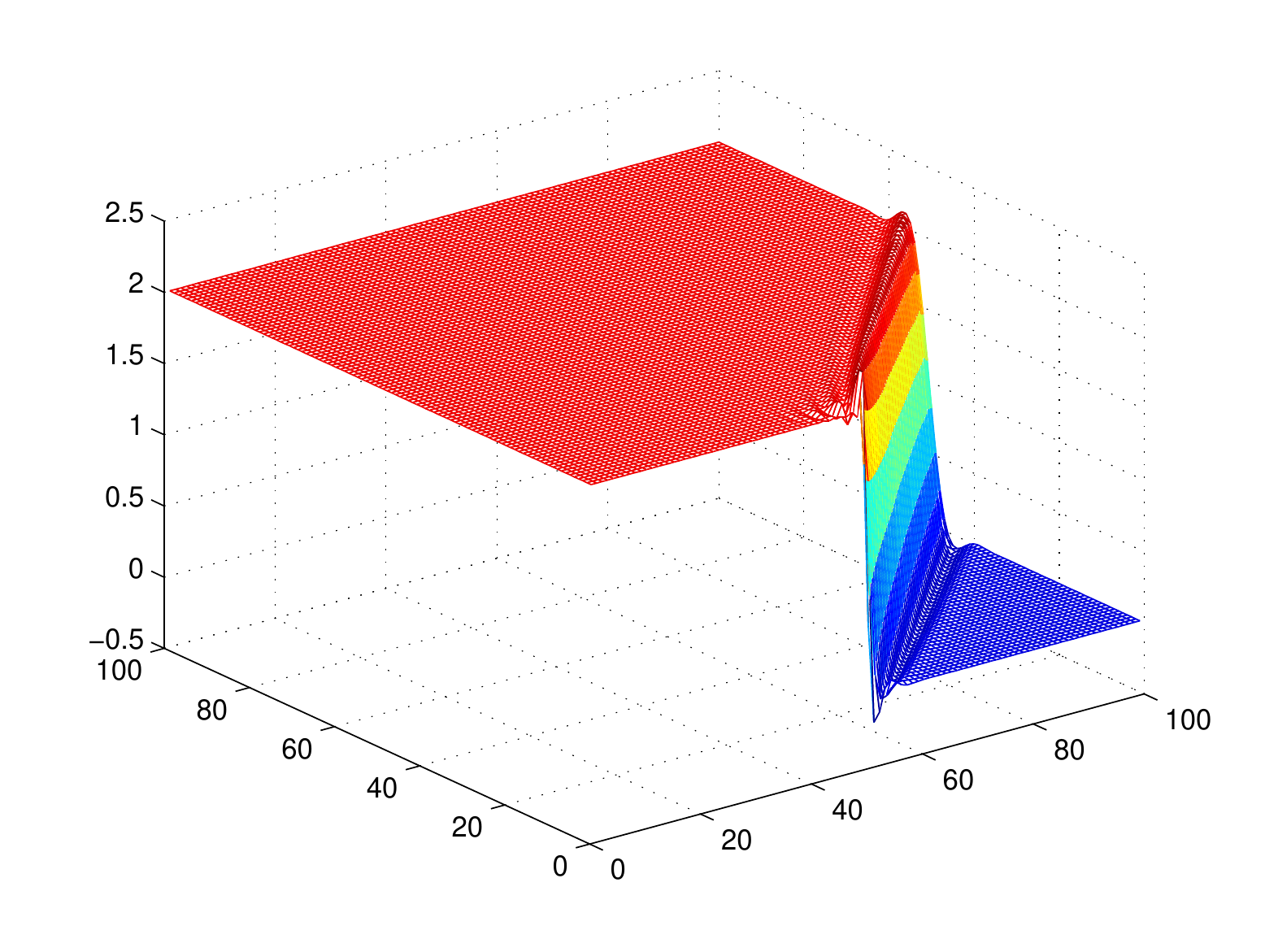}}
\subfigure[$SBP4$, fourth-order diffusion]{\includegraphics[width=0.45\linewidth]{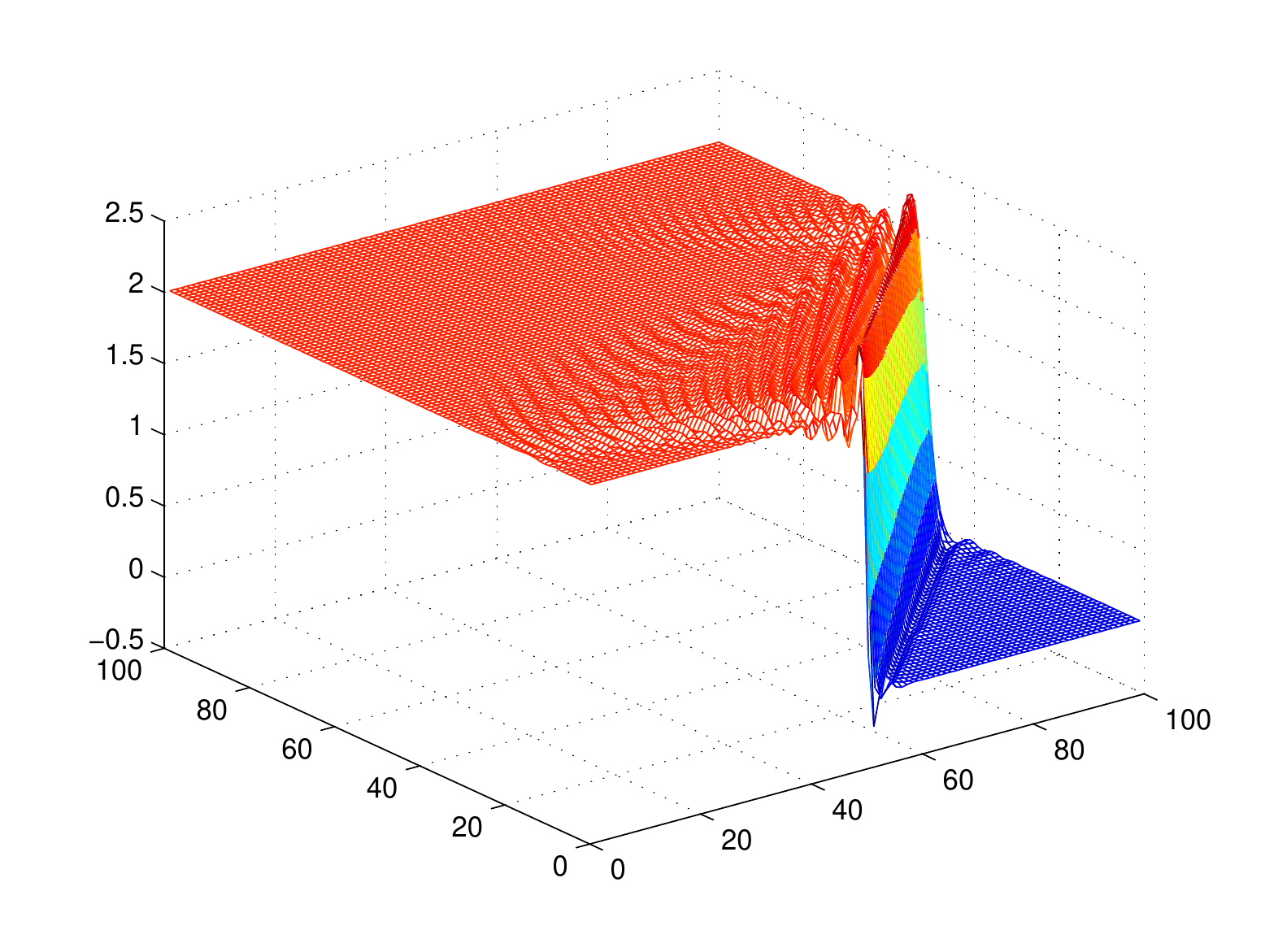}}

  \caption{Numerical results for $ B^1 (x,y,0.5)$ in experiment $3$.}
  \label{fig:5}
\end{figure}

We test the $SBP2$ ($SBP4$) scheme  with the standard second-order (fourth-order)
numerical diffusion operator as well as the scaled numerical
diffusion operator to obtain the first-order (third-order) $SBP1$ and $SBP3$ schemes.  The results on a $100 \times 100$
mesh at time $t = 0.5$ are plotted in Figure~\ref{fig:5}.
A plot at this time is of interest as some part of the solution has
interacted with the boundary and exited the domain, whereas most of the
front is still inside the domain. From Figure~\ref{fig:5}, we see that the boundary discretization works
well in all cases and does not lead to any significant oscillations in
the domain. The $SBP1$ scheme is the most dissipative with significant
smearing at the discontinuity. However, this scheme also has no
over/under shoots or oscillations and the solution is $TVD$. The $SBP2$
scheme with second-order numerical diffusion operator is oscillatory
near the discontinuity with dispersive waves on both sides of it.
The smearing is considerably less than that of the $SBP1$ scheme. The $SBP4$
scheme with standard fourth-order numerical diffusion is even more
oscillatory and leads to a larger overshoot. The $SBP3$ scheme
damps these oscillations
somewhat and still keeps the sharpness at the discontinuity making it
an acceptable alternative. 

\section{Conclusion}
\label{sec:conc}
We have considered the magnetic induction equations that arise as a
submodel in the MHD equations of plasma physics. Various forms of
these equations were presented including the symmetric forms that are
well-posed with general initial data and  Dirichlet boundary
conditions.

Standard numerical methods of the finite difference/finite volume type
have dealt with discretizations of the constrained form
(\ref{eq:Maxwell3D}) and attempted to preserve some form of the
divergence constraint. 

We describe $SBP$-$SAT$ based finite difference schemes for the initial-
boundary value problem corresponding to the magnetic induction
equations. These schemes were based on the non-conservative symmetric
form (\ref{eq:induc1}) and use $SBP$ finite difference operators to
approximate spatial derivatives and a $SAT$ technique for implementing
boundary conditions. The resulting schemes were energy stable and
higher order accurate. 

These schemes were tested on a series of numerical experiments, which
illustrated their stability and high-order of accuracy. Interesting solution features were resolved very well.
The fourth-order scheme was found to be well suited for long time
integration problems. Despite the fact that the schemes were not
preserving any particular form of discrete divergence as well as the
lack of a rigorous discrete divergence bound, the divergence errors
generated by the schemes were quite low and converged to zero at the
expected rates when the mesh was refined. The schemes were compared with
two existing lower order schemes and one divergence preserving second
order scheme. Despite lacking any divergence bounds, the $SBP$ schemes
performed at least as well as the schemes with a divergence bound. 

The numerical experiments indicate that the $SBP$-$SAT$ framework is
effective in approximating solutions of the magnetic induction
equations to a high 
order of accuracy. In the future we plan to extend these schemes to 
magnetic induction equations with resistivity.

\end{document}